\documentclass[12pt]{article}
\usepackage{amsmath}
\usepackage{amssymb}
\usepackage{amsthm}
\usepackage{amscd}
\usepackage{amsxtra}
\usepackage{verbatim}
\usepackage{xcolor}
\usepackage{color}
\usepackage{enumerate}
\usepackage{mathrsfs}
\usepackage[all]{xy} 

\usepackage{array,arydshln}
\usepackage{bm}

\allowdisplaybreaks

\numberwithin{equation}{section}

\definecolor{dblue}{rgb}{0,0,0.45}
\definecolor{red}{rgb}{0.7,0,0}

 \RequirePackage{geometry}
\newtheorem{theorem}{Theorem}[section]
\newtheorem{lemma}[theorem]{Lemma}
\newtheorem*{lemma*}{Lemma}
\newtheorem{corollary}[theorem]{Corollary}
\newtheorem{proposition}[theorem]{Proposition}

\theoremstyle{definition}

\newtheorem{remark}[theorem]{Remark}
\newtheorem{definition}[theorem]{Definition}
\newtheorem{example}[theorem]{Example}

\theoremstyle{remark}

\newcommand{\E}{{\mathbb E}}
\newcommand{\N}{{\mathbb N}}
\newcommand{\R}{{\mathbb R}}

\newcommand{\cC}{{\mathcal C}}

\newcommand{\cF}{{\mathcal F}}

\newcommand{\cH}{{\mathcal H}}

\newcommand{\cK}{{\mathcal K}}

\newcommand{\cP}{{\mathcal P}}

\newcommand{\cQ}{{\mathcal Q}}

\newcommand{\cW}{{\mathcal W}}
\newcommand{\cX}{{\mathcal X}}
\newcommand{\cY}{{\mathcal Y}}

\newcommand{\ve}{\varepsilon}

\newcommand{\la}{\langle}
\newcommand{\ra}{\rangle}
\newcommand{\nn}{\nonumber}

\newcommand{\vertiii}[1]{{\left\vert\kern-0.25ex\left\vert\kern-0.25ex\left\vert #1 
    \right\vert\kern-0.25ex\right\vert\kern-0.25ex\right\vert}}

\date{}
\begin{document}

\title{
Wong-Zakai approximation of density functions
}
\author{   
Yuzuru \textsc{Inahama} 
}
\maketitle

\begin{abstract} 
In this paper we prove the Wong-Zakai approximation
of  probability density functions of solutions at a fixed time
of rough differential equations driven by fractional 
Brownian rough path with Hurst parameter $H \in (1/4, 1/2]$.
Besides rough path theory, we use Hu-Watanabe's 
approximation theorem in the framework of 
Watanabe's distributional Malliavin calculus.
When $H=1/2$, the random rough differential equations
coincide with the corresponding
 Stratonovich-type stochastic differential equations.
Even in that case, our main result seems new.
\vskip 0.08in
\noindent {\bf Running head.} Wong-Zakai approximation of density functions.
\end{abstract}

\section{Introduction}
Consider the following stochastic differential equation (SDE)
of Stratonovich type on $\R^e$ 
driven by a standard $d$-dimensional Brownian motion 
$w =(w_t)_{t \in [0,1]}$:
\begin{equation}
dy_t = \sigma (y_t) \circ dw_t + b  (y_t)  dt, 
\qquad
\quad  y_0 =a.
\label{def.SDE}
\end{equation} 
Here, $a \in \R^e$ is an arbitrary (deterministic) starting point, 
the time interval is $[0,1]$
and the coefficients
$\sigma \colon \R^e \to \R^{e \times d}$ 
and $b \colon \R^e \to \R^e$ are assumed to be sufficiently nice.

For $m \ge 1$, denote by $w(m)=(w(m)_t)_{t \in [0,1]}$
the piecewise linear approximation of $w$ associated with 
the equal partition $\{ j/m \colon 0\le j \le m\}$ of length $1/m$  of $[0,1]$.
Let $y(m)$ be a unique solution of the corresponding 
Riemann-Stieltjes ordinary differential equation (ODE) 
driven by $w(m)$:
\begin{equation}
dy(m)_t = \sigma (y(m)_t)  dw(m)_t + b  (y(m)_t)  dt, 
\qquad
\quad  y(m)_0 =a.
\label{def.approxODE}
\end{equation} 
The famous Wong-Zakai approximation theorem 
states that 
$y(m)\to y$ as $m\to\infty$ at the process level 
in $L^q$-norm ($1<q<\infty$) with rate $m^{- (1/2 -\ve)}$
for every sufficiently small $\ve >0$, that is, 
\[
\E [ \| y(m) -y \|_\infty^q ]^{1/q} = 
O \Bigl( \Bigl( \frac{1}{m} \Bigr)^{\frac12 -\ve} \Bigr)
\qquad 
\mbox{as $m\to\infty$.}
\]
Here, $\| \cdot\|_{\infty}$ is the usual sup-norm
and $O$ is the big Landau symbol.
See \cite[Section VI-7]{iwbk} or \cite[Section 11.4]{hu} for example.
(Though it is not discussed in this paper, 
logarithmic sharpening of the convergence rate is also known.
See \cite[Theorem 11.6]{hu} for example.)

When the law of $y_t$ (for a fixed time $t$) admits a density function
$p_t(\xi)$
with respect to the Lebesgue measure $d\xi$ on $\R^e$,
can we approximate it by (functionals of) $\{y(m)\}_{m\ge 1}$?
Even though this kind of question is quite natural
and the Wong-Zakai approximation has been extensively studied, 
there seems to be no results on the above question.
In this paper we provide an answer to this question by 
combining Malliavin calculus and rough path theory (RP theory).

Since we use RP theory, the driving process 
need not be Brownian motion or its RP lift, 
but can be fractional Brownian rough path (fBRP),
that is, a canonical RP lift of 
fractional Brownian motion (fBM)
with Hurst parameter $H \in (1/4, 1/2]$.
(We denote fBM by $w = (w_t)_{t\in [0,1]}$ again.
The piecewise linear approximation $w(m)$ is defined in the same way.)
Take any $p \in (H^{-1}, [H^{-1}]+1)$ 
and let $\mathbf{w}=(\mathbf{w}^1, \ldots, \mathbf{w}^{[p]})$
be fBRP (i.e. a canonical lift of 
the $d$-dimensional fBM $w$).
We study the following rough differential equation (RDE)
driven by fBRP $\mathbf{w}$ in the $p$-variational setting:
\begin{equation}
dy_t = \sigma (y_t) dw_t + b  (y_t)  dt, 
\qquad
\quad  y_0 =a \, \in \R^e.
\label{def.introRDE}
\end{equation} 
Here, 
$\sigma \colon \R^e \to \R^{e \times d}$ 
and $b \colon \R^e \to \R^e$ are assumed to be of 
$C^\infty_{\mathrm{b}}$. 
We denote 
a unique solution of \eqref{def.introRDE} by 
$\mathbf{y}=(\mathbf{y}^1, \ldots, \mathbf{y}^{[p]})$ and  set 
$y_t = a + \mathbf{y}^1_{0,t}$.
When $H=1/2$, this process $(y_t)$ coincides with the solution of 
Stratonovich SDE \eqref{def.SDE}.
By replacing the driving process of ODE \eqref{def.approxODE}
by the piecewise linear approximation of fBM, 
we obtain an approximating process $y(m)$ in the same way.

The Wong-Zakai approximation for this random RDE 
with the reasonable convergence rate $2H-(1/2) -\ve$ was first obtained 
in \cite{fr}, but the convergence in \cite{fr} was in the pathwise sense.
The $L^q$-convergence with that rate
(i) in the case $H\in (1/3, 1/2]$ was shown in \cite{na} and 
(ii) in the case $H\in (1/4, 1/2]$ was shown in \cite{bfrs}.
(In fact,  \cite{bfrs} deals with quite general Gaussian RPs 
including fBRP with $H\in (1/4, 1/2]$.)
In this paper, however, we need  the convergence in 
every Sobolev norm to prove our main theorem
(see Theorem \ref{thm.WZdinfty} below).

Now we provide a simplified version of our main result.
(A full statement will be given in Theorem \ref{thm.main} below.)
We believe that our convergence rate is satisfactory.
First, we set 
\begin{equation}\label{def.HK}
\varphi_\rho (x) = (2\pi \rho^2)^{-e/2} e^{-|x|^2 /2\rho^2},
\qquad 
 (x, \rho) \in \R^e \times (0,\infty).
\end{equation}
This is the density function of 
the mean-zero normal distribution on $\R^e$ with 
covariance matrix $\rho^2 {\rm Id}_e$.
It should be recalled that 
if $y_t \in \mathbf{D}_\infty (\R^e)$ is non-degenerate in the sense of Malliavin, then the law of $y_t$ admits a smooth density, 
which is denoted by $p_t (\xi)$.

\begin{theorem} \label{simple}
Let the situation and the notation be as above. 
Let $t \in (0,1]$ and 
assume that $y_t \in \mathbf{D}_\infty (\R^e)$ is non-degenerate in the sense of Malliavin.
Then, for every $\delta >0$ and sufficiently small $\ve >0$, 
we have
\begin{equation}\label{eq.huwat2}
\sup_{\xi\in\R^e}  \bigl|  
  \E [ \varphi_{m^{-\delta}}  (y(m)_t - \xi)]
      - p_t (\xi)       \bigr| 
      =               
                    O \Bigl( \Bigl( \frac{1}{m} \Bigr)^{(2H- \frac12 -\ve) \wedge \delta} \Bigr)
\qquad 
\mbox{as $m\to\infty$.}                    
\end{equation}  
\end{theorem}

\begin{remark} 
As for the starting point $a$ of RDE \eqref{def.introRDE},
we may and will assume without loss of generality that $a=0$ (except in Section \ref{sec.Jac}).
\end{remark}

For high order It\^o-Taylor (or Euler) approximations of usual SDEs, 
approximations of density functions of this kind
 have been intensively studied.
For example,  see \cite{hw, bt, kh, kh01} among many others.
(In \cite{br} approximations of density functions
for various approximation schemes are proved, but 
the Wong-Zakai scheme does not seem to be included.)

Among these works, 
the one most relevant to ours is Hu-Watanabe \cite{hw}, which is well-summarized in \cite[Chapter 11]{hu}.
In that paper, they developed a general theory for 
approximations of density functions
by using Watanabe's distributional Malliavin calculus.
Simply put, if a sequence $\{F_n\}_{n\ge 1}$ 
of Wiener functionals converges to a non-degenerate
Wiener functional $F$ at a suitable rate 
in the space of test Wiener functionals,
then we automatically have an approximation theorem 
for the density of the law of $F$ (see Theorem \ref{thm.hw2.1} below
for a precise statement).

Thanks to this theorem, we have only to prove 
that, at a fixed time $t$, the Wong-Zakai approximation
holds at a suitable convergence rate 
in the space of test Wiener functionals.
This kind of convergence was already shown in \cite{ina14}
in a proof of Malliavin differentiability of RDE solutions, 
 but it was merely qualititative.
Therefore, we must make it quantitative to obtain 
an appropriate convergence rate. (See Theorem \ref{thm.WZdinfty} below,
which can be in a sense viewed as our true main result.) 
To carry it out, we will combine the following two facts:
{\rm (i)}~Convergence rate of the lift of $w(m)$ to $\mathbf{w}$ as 
$m\to\infty$ with respect to the $p$-variation topology for 
sufficiently large $p$.
This was obtained by Friz-Riedel \cite{fr}.
{\rm (ii)}~In Malliavin calculus of RDEs, 
analysis of the Jacobian RDE 
associated with a given RDE is always very important.
It is known that Lyons' continuity theorem still holds 
for the Jacobian RDE, that is, the Lyons-It\^o map of the Jacobian RDE
is locally Lipschitz continuous.
We need a suitable quantitative estimate of this
 local Lipschitz continuity. 
To the author's knowledge, however, 
no such estimates seems to be explicitly known
(although it is probably within expectation of experts). 
We will obtain one in Subsection  \ref{subsec.43}.


The organization of this paper is as follows.
In Section 2, we recall basic facts of Malliavin calculus
which will be used in this paper.
The most important one is Hu-Watanabe's approximation 
theorem for density functions (Theorem \ref{thm.hw2.1}).
In Section 3, we gather fundamental results on
geometric RPs.
Many useful estimates for RDEs are found in Section 4.
Most of them are known, but a quantitative estimate 
of local Lipschitz continuity of the Lyons-It\^o map of  
a Jacobian RDE (Proposition \ref{pr.diff_JK}) seems somewhat new
and is crucial for our purpose.
In Section 5, two important results on fBRP are recalled.
One is Friz-Riedel's convergence rate of piecewise linear approximations. 
The other is Cass-Litterer-Lyons' exponential integrability of 
``$N$-functional" associated with fBRP.
(See \cite{fr, cll}.)
The latter guarantees that the solution of the Jacobian RDE
has moments of all order.
In Section 6, 
we consider ODEs driven by piecewise linear approximations of fBM.
Mimicking arguments in \cite{ina14}, we calculate 
Malliavin derivatives of the solutions.
Combining the results obtained in previous sections, 
we show in Section 7 in a similar way to \cite{ina14}
 that the Wong-Zakai approximation holds 
at a desired convergence rate (Theorem \ref{thm.WZdinfty}),
from which our main theorem (Theorem \ref{thm.main}) immediately 
follows.

\medskip

\noindent
{\bf Notation:}~In this paper 
we will use the following notation (unless otherwise specified).
We write $\N =\{1,2, \ldots\}$.
The time interval of (rough) paths 
and stochastic processes is $[0,1]$.
We always write $\lambda_t =t$ for $0\le t \le 1$,
which is one of the simplest $\R$-valued paths.
All the vector spaces are over $\R$.

Now we will introduce the notation for some Banach spaces.
(Below, $d, e \in \N$ and $\nabla$ is the standard gradient on 
a Euclidean space.)
\begin{itemize} 
\item
The set of all continuous path $\varphi\colon [0,1] \to\R^d$ 
is denoted by $\cC (\R^d)$. 
Equipped with the usual sup-norm $\|\varphi\|_{\infty} := \sup_{0\le t \le 1}|\varphi_t|$,
this is a Banach space.
The difference of $\varphi$ is often denoted by $\varphi^1$,
that is, $\varphi^1_{s,t} := \varphi_t - \varphi_s$ for $s\le t$.
For a given starting point $a \in \R^d$, we write
$\cC_a (\R^d) =\{ \varphi \in \cC (\R^d) \colon \varphi_0 =a\}$.

\item
Let $1 \le p <\infty$. The $p$-variation seminorm of 
$ \varphi \in \cC (\R^d)$ is defined as usual by
\[
\|  \varphi\|_{p\textrm{-}\mathrm{var}} :=
\Bigl(
\sup_{\cP} \sum_{i=1}^N | \varphi^1_{t_{i-1}, t_{i}}|^p
\Bigr)^{1/p},
\qquad
\mbox{where }  \cP =\{ 0 =t_0 <t_1 <\cdots <t_N= 1\}.
\]
In the supremum above,  $\cP$ runs over all (finite) partition of $[0,1]$.
The set of all continuous paths with finite $p$-variation is 
denoted by 
\[
\cC^{p\textrm{-}\mathrm{var}} (\R^d) =\{ \varphi \in \cC (\R^d)
\colon  \|  \varphi\|_{p\textrm{-}\mathrm{var}} <\infty\}, 
\]
which is a (non-separable)
Banach space
with the norm $\|  \varphi\|_{p\textrm{-}\mathrm{var}}  +|\varphi_0|$.
For $a\in \R^d$, 
$\cC_a^{p\textrm{-}\mathrm{var}} (\R^d)$ 
is defined in an analogous way as above.

\item
Let $U$ be an open subset of $\R^m$.
For $k \in \N \cup \{0\}$,  $C^k (U, \R^n)$ denotes the set of 
$C^k$-functions from $U$ to $\R^n$.
(When $k=0$, we simply write $C (U, \R^n)$ 
instead of $C^0 (U, \R^n)$.)
The set of bounded $C^k$-functions $f \colon U\to \R^n$
whose derivatives up to order $k$ are all bounded 
is denoted by $C_{\mathrm{b}}^k (U, \R^n)$. This is a Banach space with 
$\| f\|_{C_{\mathrm{b}}^k } := \sum_{i=0}^k \|\nabla^i f\|_{\infty}$.
(Here, $ \|\cdot\|_{\infty}$ stands for the usual sup-norm on $U$.)
As usual, we set $C^\infty (U, \R^n) := \cap_{k =0}^\infty C^k (U, \R^n)$
and $C^\infty_{\mathrm{b}} (U, \R^n) := \cap_{k=0}^\infty C^k_{\mathrm{b}} (U, \R^n)$.

\item
The set of all bounded linear maps from 
a Banach space $\cX$ to another Banach space $\cY$ 
is denoted by $\mathrm{L} (\cX, \cY)$.
Equipped with the usual operator norm, 
$\mathrm{L} (\cX, \cY)$ is a Banach space.
Likewise, if $\cX_1, \ldots, \cX_n$ ($n\ge 2$) are Banach spaces,  
$\mathrm{L}^{(n)} (\cX_1, \ldots, \cX_n;\cY)$ stands for 
the Banach space of  all bounded multilinear maps
from $\cX_1\times \ldots \times \cX_n$ to $\cY$.
As usual the norm of $A \in \mathrm{L}^{(n)} (\cX_1, \ldots, \cX_n;\cY)$
is defined by 
\[
\| A\|_{\textrm{op}} := \sup \{ \| A \la v_1, \ldots, v_n \ra\|_{\cY} : 
v_i \in \cX_i, \|v_i\|_{\cX_i} \le 1 \,\, (1\le i \le n)\}.
\]
The Banach subspace of all symmetric bounded $n$-linear maps
is denoted by $\mathrm{L}^{(n)}_{\textrm{sym}} (\cX, \ldots, \cX;\cY)$.
Note that, for $A, B \in \mathrm{L}^{(n)}_{\mathrm{sym}} (\cX, \ldots, \cX ;\cY)$, $A=B$ holds if and only if 
$A \la v, \ldots, v\ra =B \la v, \ldots, v\ra$ holds for all $v \in \cX$.

\item
For brevity, we often write 
$\R^{e \times d} := \mathrm{L} (\R^d, \R^e)$
 for the set of real $e \times d$-matrices. 
 The identity matrix of size $e$ is denoted by $\mathrm{Id}_e$ or 
 simply $\mathrm{Id}$.
In this paper we  
 equip this space with the Hilbert-Schmidt norm (instead of the operator norm).
 \end{itemize}

\section{Preliminaries from Malliavin calculus}

In this section,  $(\cW, \cH, \mu)$ is 
an abstract Wiener space. 
That is, ~$(\cW, \|\cdot\|_{\cW})$~is a separable Banach space, ~$(\cH, \|\cdot\|_{\cH})$~ is a separable Hilbert space, $\cH$~is a dense subspace of~$\cW$~and the inclusion map is continuous, and~$\mu$~is the (necessarily unique) probability measure on~$(\cW, \mathcal{B}_{\cW})$~with the property that
\begin{equation}\label{abspro}
\int_{\cW}\exp\Bigl(\sqrt{-1}_{\cW^*}\langle \psi, w \rangle_{\cW}\Bigr)\mu (d w)=\exp\Bigl(-\frac{1}{2}\|\psi\|^2_{\cH^*}\Bigr),
\qquad
\qquad
\psi \in \cW^* \subset \cH^*,
\end{equation}
where we have used the fact that $\cW^*$ becomes a dense subspace of $\cH$ when we make the natural identification between $\cH^*$ and $\cH$ itself. 
Hence, $\cW^* \hookrightarrow \cH^* =\cH \hookrightarrow\cW$
and both inclusions are continuous and dense.
(For basic information on abstract Wiener spaces, see \cite{sh, hu, str} among others.)


We first set the notation and 
quickly summarize some basic facts in Malliavin calculus
which are related to Watanabe distributions
(i.e. generalized Wiener functionals).
Most of the contents and the notation
in this section are found in 
 \cite[Sections V.8--V.10]{iwbk} with trivial modifications.
Also, \cite{sh, nu, hu, mt} are good textbooks on Malliavin calculus.

As for differential operators on $(\cW, \cH, \mu)$,
$D$ stands for the $\cH$-derivative (i.e. the gradient operator 
in the sense of Malliavin calculus), while 
$L = -D^* D$ stands for the Ornstein-Uhlenbeck operator.

We denote by $\mathscr{S}({\mathbb R}^e)$
and $\mathscr{S}^{\prime}({\mathbb R}^e)$ 
the Schwartz class of smooth rapidly decreasing functions 
and its dual (i.e. the set of tempered Schwartz distributions), respectively.
Below, $e \in \N$ 
and ${\cal K}$ is a real separable Hilbert space. 
We denote by ${\rm Pol} (\cK)$ the set of $\cK$-valued polynomials 
on $\cW$.
(For a precise definition of ${\rm Pol} (\cK)$, see \cite[Section 2]{sh}.)

\begin{enumerate}
\item[{\bf (a)}]
Sobolev spaces ${\bf D}_{q,r} ({\cal K})$ of ${\cal K}$-valued 
(generalized) Wiener functionals,
where $q \in (1, \infty)$ and $r \in {\mathbb R}$.
This space is defined as the closure of ${\rm Pol} (\cK)$
with respect to the norm 
\[
\| G \|_{{\bf D}_{q,r}} := \| (I -L)^{r/2} G \|_{L^q},
\qquad G \in {\rm Pol} (\cK).
\]
As usual, we will use the spaces 
${\bf D}_{\infty} ({\cal K})= \cap_{k=1 }^{\infty} \cap_{1<q<\infty} {\bf D}_{q,k} ({\cal K})$, 
$\tilde{{\bf D}}_{\infty} ({\cal K}) 
= \cap_{k=1 }^{\infty} \cup_{1<q<\infty}  {\bf D}_{q,k} ({\cal K})$ of test functionals 
and  the spaces ${\bf D}_{-\infty} ({\cal K}) = \cup_{k=1 }^{\infty} \cup_{1<q<\infty} {\bf D}_{q,-k} ({\cal K})$, 
$\tilde{{\bf D}}_{-\infty} ({\cal K}) = \cup_{k=1 }^{\infty} \cap_{1<q<\infty} {\bf D}_{q,-k} ({\cal K})$ of 
 Watanabe distributions as in \cite{iwbk}.
When ${\cal K} ={\mathbb R}$, we simply write ${\bf D}_{q, r}$, etc.
Recall that the Sobolev norms (resp. spaces)
are non-decreasing (resp. non-increasing) in both $q$ and $r$.
\item[{\bf (b)}]  
Meyer's equivalence of Sobolev norms.
(See \cite[Theorem 8.4]{iwbk}. 
A stronger version can be found in \cite[Theorem 4.6]{sh} or \cite[Theoem 1.5.1]{nu}.)  It states that,
   for every $q \in (1, \infty)$ and $k \in {\mathbb N}$, there exists 
   a  constant $C= C_{q,k}\ge 1$ such that
   \[
   C^{-1} \|G\|_{{\bf D}_{q,k}}  \le  \|G\|_{L^q} +  \|D^kG\|_{L^q}   \le C\| G \|_{{\bf D}_{q,k}} ,
\qquad G \in {\rm Pol} (\cK).
\]
\item[{\bf (c)}] 
For 
$F =(F^1, \ldots, F^e) \in {\bf D}_{\infty} ({\mathbb R}^e)$, we denote by 
$\Sigma^{ij}_F (w) =  \la DF^i (w),DF^j (w)\ra_{{\cal H}}$
 the $(i,j)$-component of Malliavin covariance 
 matrix ($1 \le i,j \le e$).
We say that $F$ is non-degenerate in the sense of Malliavin
if $(\det \Sigma_F)^{-1} \in  \cap_{1<q< \infty} L^q (\mu)$.
If $F \in {\bf D}_{\infty} ({\mathbb R}^e)$
is non-degenerate, its law on $\R^e$ admits a smooth
rapidly decreasing
density function $p_F =p_F (y)$ with respect to the Lebesgue measure $dy$,
that is, $(\mu \circ F^{-1} )(dy)= p_F (y)dy$.
\item
[{\bf (d)}] Pullback $T \circ F =T(F)\in \tilde{\bf D}_{-\infty}$ of
 $T \in \mathscr{S}^{\prime}({\mathbb R}^e)$
by a non-degenerate Wiener functional $F \in {\bf D}_{\infty} ({\mathbb R}^e)$. 
The most important example of $T$ is Dirac's delta function.
In that case, ${\mathbb E}[\delta_y (F)] :=\la \delta_y (F), \mathbf{1}\ra =p_F (y)$ 
holds for every $y\in\R^e$.
Here, $\la \star, *\ra$ denotes the pairing of ${\bf D}_{-\infty}$ and ${\bf D}_{\infty}$ as usual and $\mathbf{1}$ is the constant function $1$.
(See \cite[Section 5.9]{iwbk}.)
\end{enumerate}

\begin{remark} 
In some of the books cited above (in particular \cite{iwbk, mt}),  results are formulated on  a special Gaussian space.
However, almost most all of them 
(at least, those that will be used in this paper) still hold true
 on any abstract Wiener space.
\end{remark}


Now let us recall a few important properties of the Wiener chaos.
For $n\in \N\cup\{0\}$, denote by
$\mathscr{C}_n$ the $n$th homogeneous Wiener chaos
(see \cite[Chapter 1]{sh} for example).
It is well-known that $\mathscr{C}_0$ is the space of 
constant functions
and $\mathscr{C}_1 = \{ \langle k, \bullet \rangle \colon k\in \cH\}$
i.e. the family of centered Gaussian random variables 
defined on $\cW$ indexed by $\cH$.
(It should be noted that $\langle k, \bullet \rangle$ is defined as 
a Wiener integral if it does not belong to $\cW^*$)
The $n$th inhomogeneous Wiener chaos is denoted by
 $\mathscr{C}^{\prime}_n := \oplus_{i=0}^n  \mathscr{C}_i$.
Heuristically, $\mathscr{C}^{\prime}_n$ is the set of 
all real-valued
``polynomials" of order at most $n$ on $\cW$.
The most important fact is the orthogonal decomposition
$L^2 (\mu) =  \oplus_{i=0}^\infty  \mathscr{C}_i$.
Moreover, $\mathscr{C}_n$ is the eigenspace of $-L$ in $L^2 (\mu)$ associated 
with the eigenvalue $n$.
If each component of $\R^d$-valued function $F =(F^1, \ldots, F^d)$ on $\cW$ 
belongs to $\mathscr{C}_n$ (resp. $\mathscr{C}^\prime_n$),
we say that $F$ belongs to $\mathscr{C}_n (\R^d)$
(resp. $\mathscr{C}^\prime_n (\R^d)$).

Restricted to a fixed homogeneous Wiener chaos $\mathscr{C}_n$,
all the $L^q$-norms ($1 <q <\infty$) are known to be equivalent,
that is,
there exists a constant $C=C_{n,q} \ge 1$ such that
\begin{equation}\label{eq.0111-1}
C^{-1} \| F\|_{L^2}  \le  \| F\|_{L^q} \le  C \| F\|_{L^2},
\qquad
F \in \mathscr{C}_n.
\end{equation}
This is a consequence of the hypercontractivity.
For details, see \cite[Theorem 2.14]{sh} for example.
%
%
In this paper we use this equivalence in the following form.
\begin{lemma} \label{lem.equiv_chaos}
Let $n\in \N \cup\{0\}$, $2 \le q <\infty$ and $k \in \N\cup\{0\}$.
Then, restricted to $\mathscr{C}^\prime_n$, 
all $\mathbf{D}_{q,k}$-norms are equivalent, that is, 
there exists a constant $C=C_{n,q, k} >0$ satisfying that
\[
\| F\|_{L^2}  \le  \| F\|_{\mathbf{D}_{q,k}} \le  C \| F\|_{L^2},
\qquad
F \in \mathscr{C}^\prime_n.
\]
\end{lemma}

\begin{proof} 
The left inequality is trivial. We show the right one.
In this proof, the positive constant $C$ may change from line to line.

Denote by $\Pi_n \colon L^2 (\mu) \to \mathscr{C}_n$
the orthogonal projection and write $F = \sum_{i=0}^n \Pi_i F$ 
for $F \in \mathscr{C}^\prime_n$. 
Then, we have
\begin{align*}
\| F \|_{\mathbf{D}_{q,k}} 
  := \| (I-L)^{k/2} F \|_{L^q}
     &\le C \sum_{i=0}^n (1+i)^{k/2}  \| \Pi_i  F \|_{L^q}
         \\
         &\le C \sum_{i=0}^n (1+i)^{k/2}  \| \Pi_i  F \|_{L^2}
          \\
            &\le C (1+n)^{k/2}  \sqrt{n}
               \Bigl( \sum_{i=0}^n \| \Pi_i  F \|_{L^2}^2  \Bigr)^{1/2} 
                  = C  \|F \|_{L^2}.
                                          \end{align*}
Here, we used \eqref{eq.0111-1} and
the fact that each $F \in \mathscr{C}_n$ is an eigenfunction of 
$L$ with eigenvalue $n$.  
\end{proof}

\begin{lemma} \label{lem.L_sym}
Let $n \in \N$.
Then, there exists a constant $C=C_{n}>0$ satisfying that
\[
\|A\|_{\cH^{* \otimes n}}^2
  \le C \int_{\cW}
       |A\la w, \ldots, w \ra |^2   \mu(dw), \qquad
         A \in  \mathrm{L}^{(n)}_{\mathrm{sym}} (\cW, \ldots, \cW ;\R).
  \]
\end{lemma}

\begin{proof} 
Set $F (w) =  A\la w, \ldots, w \ra$. Then, one can easily see that
$F\in \mathscr{C}^\prime_n$ and $D^n F (w) =n! A$ (constant in $w$).
Since $ n!\|A\|_{\cH^{* \otimes n}} =  
\E[ \|D^n F\|_{\cH^{* \otimes n}}^2 ]^{1/2} \le \|F\|_{ \mathbf{D}_{2,n}}$,
the desired inequality immediately follows from Lemma \ref{lem.equiv_chaos} above.
\end{proof}

\begin{remark} \label{rem.L_sym}
In Lemmas \ref{lem.equiv_chaos} and \ref{lem.L_sym}, 
$F$ and $A$ are real-valued.
However, it immediately follows that these lemmas 
still hold for the $\R^d$-valued case.
(Of course, the constant $C$ may depend on $d$.)
\end{remark}


The following theorem is \cite[Theorem 2.1]{hw}.
We denote the standard Laplacian on $\R^e$ by $\Delta$.
Recall that $\varphi_\rho (x)$ is defined by \eqref{def.HK}.

\begin{theorem} \label{thm.hw2.1}
Let $F \in {\bf D}_{\infty} ({\mathbb R}^e)$ be non-degenerate 
in the sense of Malliavin and let $F_m \in {\bf D}_{\infty} ({\mathbb R}^e)$
for $m\in \N$.
Suppose that $F_m$ approximates  $F$ in ${\bf D}_{\infty} ({\mathbb R}^e)$
  with order $\gamma >0$, that is, for every $q \in (1,\infty)$ 
    and $r >0$, 
    \[
      \| F_m -F \|_{ \mathbf{D}_{q,r}} = O (m^{-\gamma})
    \]
as $m\to\infty$.
Here, $O$ is the big Landau symbol as usual.

Then, for every $r>0$, $\beta \ge 0$, $\delta >0$ and $1<q< \infty$
  satisfying that $r > \beta + e/q' +1$ (with $1/q +1/q' =1$), it holds that
\begin{equation}\label{eq.huwat}
\sup_{y\in\R^e}  \bigl\|  
  [(1 -\Delta)^{\beta/2}\varphi_{m^{-\delta}}] (F_m -y)
      - [(1 -\Delta)^{\beta/2}\delta_y] (F)
       \bigr\|_{ \mathbf{D}_{q, -r}} = 
                    O (m^{-\gamma \wedge \delta})
\end{equation}  
as $m\to\infty$.
\end{theorem}

Note that in Theorem \ref{thm.hw2.1} above, 
the non-degeneracy of $F_m$ is not assumed.
Note also that $|\E [\Phi] | := |\langle \Phi, \mathbf{1}\rangle|
\le \| \Phi\|_{\mathbf{D}_{q, -r}}$ for all $q\in (1,\infty)$, $r \ge 0$ and 
$\Phi \in \mathbf{D}_{q, -r}$.
Hence, when $\beta =0$ and $r >0$ is large enough in the above theorem, 
${\mathbb E}[\delta_y (F)] = p_F (y)$ is approximated by 
$\{\varphi_{m^{-\delta}} (F_m -y)  \}_{m=1}^\infty$.

\section{Geometric rough paths}

In this section we summarize some known facts in RP theory
such as the RP lift, Lyons' extension theorem and 
RP integration. 
(Discussions on RDEs will be postponed to 
the next section.)
Everything in this and the next sections is deterministic.
Throughout this section, we assume $2 \le p <\infty$ and $d\in \N$.

Set $\triangle :=\{ (s,t) : 0\le s\le t \le 1  \}$.
A continuous  map $\omega \colon \triangle \to [0,\infty)$ 
is called a control function if it is superadditive, that is,
$\omega (s,u) +  \omega (u, t) \le  \omega (s, t)$ for all $s \le u \le t$.
For $\rho \in [1,\infty)$ and
a continuous map $A\colon \triangle \to \R^d$ that vanish 
identically on the diagnal, we set 
\[
\|  \varphi\|_{\rho\textrm{-}\mathrm{var}} :=
\Bigl(
\sup_{\cP} \sum_{i=1}^N | A_{t_{i-1}, t_{i}}|_{\R^d}^\rho
\Bigr)^{1/\rho},
\qquad 
\mbox{where }
\cP =\{ 0 =t_0 <t_1 <\cdots <t_N= 1\}.
\]
In the supremum above,  $\cP$ runs over all (finite) partition of $[0,1]$.
For $n \in \N$,
the truncated tensor algebra over $\R^d$ of degree $n$
is denoted by $T^{(n)}(\R^d) = \oplus_{i=0}^n (\R^d)^{\otimes i}$,
where $(\R^d)^{\otimes 0} :=\R$.

Let $p \in [2,\infty)$.
A continuous map 
$\mathbf{x}= (1, \mathbf{x}^1, \ldots, \mathbf{x}^{[p]})
\colon \triangle \to T^{([p])}(\R^d)$
is called an RP over $\R^d$ (of roughness $p$) 
if the following two conditions are 
satisfied:
\begin{align} 
\mathbf{x}_{s,t} = \mathbf{x}_{s,u}\otimes \mathbf{x}_{u,t},
\qquad 
 s \le u \le t.
 \label{eq.0112-1}
 \\
 \max_{1 \le i \le [p]} \|  \mathbf{x}^i \|_{p/i\textrm{-}\mathrm{var}} <\infty.
 \label{eq.0112-2}
\end{align}
The set of all RPs of roughness $p$ is denoted by $\Omega_p (\R^d)$,
which is a complete metric space with the distance 
\[
d_{p\textrm{-}\mathrm{var}} (\mathbf{x}, \hat{\mathbf{x}}) := \max_{1 \le i \le [p]} \|  \mathbf{x}^i -  \hat{\mathbf{x}}^i \|_{p/i\textrm{-}\mathrm{var}}.
\]
Condition \eqref{eq.0112-1}, which is called Chen's identity, 
guarantees that $\mathbf{x}^i$
vanishes identically on the diagonal.
More concretely, \eqref{eq.0112-1} is equivalent to 
\[
\mathbf{x}^k_{s,t} = \sum_{i=0}^k 
\mathbf{x}^{k-i}_{s,u}\otimes \mathbf{x}^{i}_{u,t},
\qquad 
1 \le k \le [p], \,\,
 s \le u \le t.
 \]

If $x \in \cC^{\rho\textrm{-}\mathrm{var}}_0 (\R^d)$ for some $\rho \in [1,2)$,
we can define 
\begin{equation}\label{def.itr.Young}
\mathbf{x}^k_{s,t} := \int_{s \le t_1 \le \cdots \le t_k \le t} 
  dx_{t_1} \otimes \cdots \otimes dx_{t_k},
  \qquad 
     (s,t) \in  \triangle
   \end{equation}
for all $k\in \N$ as an iterated Young integral.
Obviously, $\mathbf{x}^1_{s,t} =x_t -x_s$
and $\mathbf{x}^{k+1}_{s,t} = \int_s^t \mathbf{x}^{k}_{s,u} \otimes dx_u$.
For every $p \in [2,\infty)$,
$S_p (x) :=(1, \mathbf{x}^1, \ldots, \mathbf{x}^{[p]}) \in \Omega_p (\R^d)$,
which is called the natural lift of $x$.
 For every $\rho$ and $p$ as above, the lift map 
    $S_p \colon \cC^{\rho\textrm{-}\mathrm{var}}_0 (\R^d)\to \Omega_p (\R^d)$
      a locally Lipschitz continuous injection.
      (We say a map between two metric spaces is 
       {\it locally Lipschitz continuous} if the map is Lipschitz continuous
         when restricted to every bounded subset of the domain.)

        For $p \in [2,\infty)$, 
         we define the geometric RP space by 
         \[
         G\Omega_p (\R^d) := \overline{ \{ S_p (x) \,:\,  x \in \cC^{\rho\textrm{-}\mathrm{var}}_0 (\R^d)\}},
                  \]
where the closure is taken with respect to $d_p$.
It is known that the right hand side does not depend on $\rho \in [1,2)$.
Elements of $G\Omega_p (\R^d)$ are called geometric RPs over $\R^d$ (of roughness $p$).
By the way it is defined, $G\Omega_p (\R^d)$ is complete and separable.

\begin{remark} 
In the above definitions of the variation norm and
of RP spaces, the time interval is $[0,1]$. 
Obviously, 
it can easily be replaced by any subinterval $[s,t]\subset [0,1]$.
In that case, we write $\|\cdot\|_{p\textrm{-}\mathrm{var}, [s,t]}$
and $G\Omega_p ([s,t], \R^d)$, etc. 
\end{remark}

The intrinsic control function associated with $\mathbf{x}\in \Omega_p (\R^d)$ is defined by
\[
\vertiii{\mathbf{x}}_{p\textrm{-}\mathrm{var}, [s,t]}^p :=  
    \sum_{1 \le i \le [p]} \|  \mathbf{x}^i \|_{p/i\textrm{-}\mathrm{var}, [s,t]}^{p/i},
\qquad   (s,t) \in  \triangle.
\]
Obviously, this controls $\mathbf{x}$, that is,
$|\mathbf{x}^i_{s,t}| \le (\vertiii{\mathbf{x}}_{p\textrm{-}\mathrm{var}, [s,t]}^p)^{i/p}$ for all
$(s,t) \in  \triangle$ and $1 \le i \le [p]$.
The homogeneous norm of $\mathbf{x}$ is defined by 
  $\vertiii{\mathbf{x}}_{p\textrm{-}\mathrm{var}}:= \vertiii{\mathbf{x}}_{p\textrm{-}\mathrm{var}, [0,1]}$
(which is slightly different from, but equivalent to the standard definition).


For $\mathbf{x}\in \Omega_p (\R^d)$ and $\beta >0$, 
set $\tau_0 =0$ and  $\tau_{m}$, $m\in\N$, inductively as follows.
\[
\tau_m := 1 \wedge \inf \{ t > \tau_{m-1} \,:\, 
\vertiii{\mathbf{x}}_{p\textrm{-}\mathrm{var}, [\tau_{m-1} ,t]}^p \ge \beta\}.
\]
Then, we set 
\begin{equation}\label{def_N}
N^p_{\beta} (\mathbf{x})  : = \max\{ m \,:\,\tau_m <1\}. 
\end{equation}
This quantity plays a very important role when we estimate moments of 
solutions of linear RDEs driven by Gaussian RPs.
By the superadditivity, 
$N^p_{\beta} (\mathbf{x}) \le \beta^{-1} \vertiii{\mathbf{x}}^p_{p\textrm{-}\mathrm{var}, [0,1]}$ holds.


We recall a few basic operations on $G\Omega_p (\R^d)$.
 In the itemizations below, $d, d' \in\N$,
 $p \in [2,\infty)$ and $\rho \in [1,2)$ with
  $p^{-1} + \rho^{-1}>1$.
\begin{enumerate} 
\item[(1)]
The following map 
\[ 
\R \times G\Omega_p (\R^d) \,\ni
   (c, \mathbf{x}) \mapsto c \mathbf{x} := 
   (1, c^1\mathbf{x}^1, \ldots, c^{[p]}\mathbf{x}^{[p]})
    \,\in G\Omega_p (\R^d)
    \]
     is locally Lipschitz continuous and called the dilation.
     In a similar way, for $A \in\mathrm{L} (\R^d, \R^{d'})=\R^{d' \times d}$,
     The following map 
\[ 
\mathrm{L} (\R^d, \R^{d'})\times G\Omega_p (\R^d) \,\ni
   (A, \mathbf{x}) \mapsto \Gamma_A \mathbf{x} := 
   (1, A\mathbf{x}^1, \ldots, A^{\otimes [p]}\mathbf{x}^{[p]})
    \,\in G\Omega_p (\R^{d'})
    \]
     is locally Lipschitz continuous. Here, 
     $A^{\otimes i}\colon (\R^d)^{\otimes i}
     \to (\R^{d'})^{\otimes i}$
     stands for the $i$-fold tensor product of $A$. 
      (In this paper, $\Gamma_A$ is called the generalized dilation by $A$)
\item[(2)] 
The addition map
$
 \cC^{\rho\textrm{-}\mathrm{var}}_0 (\R^d)\times   \cC^{\rho\textrm{-}\mathrm{var}}_0 (\R^d)
 \,\ni    (x, k) \mapsto x+k \, \in  \cC^{\rho\textrm{-}\mathrm{var}}_0 (\R^d)
$ 
   extends uniquely to a continuous map 
     $T\colon G\Omega_p (\R^d)\times   \cC^{\rho\textrm{-}\mathrm{var}}_0 (\R^d)
       \to G\Omega_p (\R^d)$
         satisfying that $T(S_p (x) ,k) =S_p (x+k)$ for all $x, k \in  \cC^{\rho\textrm{-}\mathrm{var}}_0 (\R^d)$.
           This map is locally Lipschitz continuous
            and called the Young translation.  
               We often write $T_k (\mathbf{x})$ for  $T (\mathbf{x}, k)$.

\item[(3)] 
Similarly, there exists a unique continuous map
$\textbf{Pair}\colon G\Omega_p (\R^d)\times   \cC^{\rho\textrm{-}\mathrm{var}}_0 (\R^{d'}) \to G\Omega_p (\R^d \oplus \R^{d'})$
satisfying that $\textbf{Pair} (S_p (x) ,k) =S_p ((x, k))$ for all $x \in  \cC^{\rho\textrm{-}\mathrm{var}}_0 (\R^d)$ and 
       $k \in  \cC^{\rho\textrm{-}\mathrm{var}}_0 (\R^{d'})$.
            This map is locally Lipschitz continuous
            and called the Young pairing.   
              For brevity, we will basically write $(\mathbf{x}, \mathbf{k})$
                (both are boldface letters)
               for $\textbf{Pair} (\mathbf{x}, k)$.
\end{enumerate}
To calculate (2) and (3) above, we need Young integration theory 
and the shuffle relations for (components of) geometric RPs. 
For example, the first and second level paths of 
$(\mathbf{x}, \mathbf{k})$ are given by 
$(\mathbf{x}, \mathbf{k})^1_{s,t} = (\mathbf{x}^1_{s,t}, \mathbf{k}^1_{s,t})$
and 
\[
(\mathbf{x}, \mathbf{k})^2_{s,t} 
 = 
  \Bigl(   \mathbf{x}^2_{s,t},  \,\, \int_s^t  \mathbf{x}^1_{s,u} dk_u, \,\,
     \int_s^t  \mathbf{k}^1_{s,u} d \mathbf{x}^1_{0,u}, \,\,
            \mathbf{k}^2_{s,t}       \Bigr).
\]
Here, $\mathbf{k}^i_{s,t}$ ($i=1,2$) is given in \eqref{def.itr.Young}
  and the two cross integrals are in the Young sense.


Now we review Lyons' extension theorem.
For a given $\mathbf{x}= (1, \mathbf{x}^1, \ldots, \mathbf{x}^{[p]})\in \Omega_p (\R^d)$, $p \in [2,\infty)$,
we can construct $\mathbf{x}^k \colon \triangle \to (\R^d)^{\otimes k}$ for all $k \ge [p] +1$ inductively so that 
$(1, \mathbf{x}^1, \ldots, \mathbf{x}^{k})$ is multiplicative
(i.e. it satisfies \eqref{eq.0112-1} with $[p]$ being replaced by $k$) 
for all $k$ as follows
 (see \cite[Theorem 3.1.2]{lq}).
 
  If $(1, \mathbf{x}^1, \ldots, \mathbf{x}^{k})$ is obtained
  for $k \ge [p]$, then we set
   \begin{equation} \label{def.LET} 
   \mathbf{x}^{k+1}_{s,t} :=
    \lim_{|\cP| \searrow 0}
      \sum_{i=1}^k  \sum_{l=1}^N   \mathbf{x}^{k+1-i}_{s,t_{l-1}}
         \otimes      \mathbf{x}^{i}_{t_{l-1},t_{l}},
               \end{equation}
where 
$\cP =\{ s =t_0 <t_1 <\cdots <t_N= t\}$ is a partition of $[s,t]$ and $|\cP|$
 is its mesh.
Moreover, the following estimate (called the neoclassical inequality) 
holds for $\mathbf{x}^{k}$ ($k \ge [p]+1$):
Suppose that $\omega$ is a control function such that 
 \begin{equation} \label{eq.0113-1} 
|\mathbf{x}^k_{s,t}| 
 \le \frac{ \omega (s,t)^{k/p}}{\beta_p \, (k/p)!},
\qquad\quad
(s,t) \in  \triangle, \,\, 1 \le k \le [p]
\end{equation}
holds, 
where $\beta_p := 1 + p^2 \{1 + \sum_{m=3}^\infty \{2/(m-2)\}^{([p]+1)/p}  \} >0$ is a positive constant.
  Then, \eqref{eq.0113-1} actually holds  for all $k \ge [p] +1$, too.
Since a constant multiple of  
the intrinsic control function associated with $\mathbf{x}$ satisfies 
\eqref{eq.0113-1}, $(1, \mathbf{x}^1, \ldots, \mathbf{x}^{[r]})\in \Omega_r (\R^d)$ for all $r \in [p, \infty)$.
      We write $\mathbf{Ext}_{p,r} ((1, \mathbf{x}^1, \ldots, \mathbf{x}^{[p]})) 
           = (1, \mathbf{x}^1, \ldots, \mathbf{x}^{[r]})$.
             As a map from $\Omega_p (\R^d)$ to $\Omega_r (\R^d)$,  
               this map is locally Lipschitz continuous (see \cite[Theorem 3.1.3]{lq}).
(By abusing the notation, we will sometimes denote $\mathbf{Ext}_{p,r}(\mathbf{x})$ by $\mathbf{x}$ again.) 
 
\begin{lemma} 
Let $2 \le p \le r \le r' <\infty$. Then, $\mathbf{Ext}_{p,r}$
maps $G\Omega_p (\R^d)$ to $G\Omega_r (\R^d)$
and is locally Lipschitz continuous.
Moreover,  as maps between geometric RP spaces, 
$\mathbf{Ext}_{p,r'} = \mathbf{Ext}_{r, r'} \circ \mathbf{Ext}_{p,r}$
holds.
\end{lemma}

\begin{proof} 
We will prove that, for
 $x \in \cC^{1\textrm{-}\mathrm{var}}_0 (\R^{d})$,
  $\mathbf{x}^k$ defined by \eqref{def.itr.Young} and 
       $\mathbf{x}^k$ defined by \eqref{def.LET} coincides. 
         (To distinguish, we denote the latter by $\hat{\mathbf{x}}^k$ in this proof.)
         Once this is done, the rest is trivial since  $\mathbf{Ext}_{p,r}$ is continuous and  
                 $S_p (\cC^{1\textrm{-}\mathrm{var}}_0 (\R^{d}) )$ is dense 
                   in $G\Omega_p (\R^d)$.

         Let $\omega$ be any control function such that \eqref{eq.0113-1} 
         holds for $k=1$. Then, the neoclassical inequality tells that 
          there is a constant $\kappa >0$ such that 
           $|\hat{\mathbf{x}}^k_{s,t}|        \lesssim \omega (s,t)^{1+\kappa}$
               for all $k \ge 2$.  
                 Since $\omega$ is uniformly continuous on $\triangle$ 
                 and vanishes on the diagonal, 
                   $\lim_{\delta \searrow 0} M(\delta, \omega)=0$,
                   where $M(\delta, \omega) :=\sup\{ \omega (s,t) : t-s \le \delta \}$.
                     This and the superadditivity imply that 
                     \[
                         \sum_{l=1}^N
                         | \hat{\mathbf{x}}^{k}_{t_{l-1},t_{l}}|
                         \lesssim 
                                         \sum_{l=1}^N\omega (t_{l-1},t_{l})^{1+\kappa}
                                                  \lesssim 
                                                      \omega (0,1) M(\delta, |\cP|)^\kappa\to 0
                                                               \qquad \mbox{as $|\cP| \to 0$}.              
                                                                                                    \]
This means that in the sum over $i$ in \eqref{def.LET}, 
the contribution from 
  the terms $i \neq 1$ is $0$. 
    Hence, in this case \eqref{def.LET} reads:
      \[
        \hat{\mathbf{x}}^{k+1}_{s,t} :=
    \lim_{|\cP| \searrow 0}
      \sum_{l=1}^N  \hat{\mathbf{x}}^{k}_{s,t_{l-1}}
         \otimes   (x_{t_{l}} - x_{t_{l-1}})
          = \int_s^t  \hat{\mathbf{x}}^{k}_{s,u} \otimes dx_u.     
           \] 
            The last equality is due to the definition of Young integration.
             Since $\mathbf{x}^{1}= \hat{\mathbf{x}}^{1}$, we can show
                $\mathbf{x}^{k}= \hat{\mathbf{x}}^{k}$
               for all $k$ by mathematical induction and the above fact.  
               Note that we have essentially shown that
                 $\mathbf{Ext}_{p,r} \circ S_p (x)= S_r (x)$ 
                   for $x \in \cC^{1\textrm{-}\mathrm{var}}_0 (\R^{d})$.
               \end{proof}


Now we recall RP integration.
Let $d, e \in \N$, $p \in [2, \infty)$ and $f \colon \R^e \to \R^{e \times d}$.
As is well-known, 
if $f$ is of $C^{[p] +1}$, then the Riemann-Stieltjes 
integration map 
\[
x \in  \cC^{1\textrm{-}\mathrm{var}}_0 (\R^{d}) 
 \quad \mapsto \quad \int_0^{\cdot} f(x_s) dx_s
\in  \cC^{1\textrm{-}\mathrm{var}}_0 (\R^{e}) 
\]
uniquely extends to a continuous map 
\begin{equation}\label{eq.0126-1}
\mathbf{x} \in G\Omega_p (\R^d)   \quad \mapsto  \quad
\int f(\mathbf{x}) d\mathbf{x} \in G\Omega_p (\R^d),
\end{equation}
which is called RP integration.
Moreover, this extended map is locally Lipschitz continuous.
(By ``extend", we mean $\int f(S_p(x)) dS_p(x) = 
S_p(\int_0^{\cdot} f(x_s) dx_s)$ for all $x$.)
One can easily show that RP integration is consistent is $p$,
that is, if $2\le p \le r <\infty$ and $f$ is of $C^{[r] +1}$, then 
$\mathbf{Ext}_{p,r} \circ I^f_p = I^f_r \circ \mathbf{Ext}_{p,r}$,
where $I^f_p \colon G\Omega_p (\R^d)\to G\Omega_p (\R^d)$ 
is the above RP integration map in the $p$-variation setting
in \eqref{eq.0126-1}.

RP integration is known to satisfy the following estimates:
\begin{proposition} \label{pr.RPint}
Let the situation be as above and write 
$\mathbf{y}=\int f(\mathbf{x}) d\mathbf{x}$ for simplicity.
\\
{\rm (i)}~If $f$ is of $C_{\mathrm{b}}^{[p] +1}$, then there exists
a constant $c>0$ (independent of $\mathbf{x}, \hat{\mathbf{x}}$ and $f$) such that
\begin{align} 
\vertiii{\mathbf{y}}_{p\textrm{-}\mathrm{var}} 
&\le  
c (1 + \| f\|_{C_{\mathrm{b}}^{[p] +1}})^c (1+ \vertiii{\mathbf{x}}_{p
\textrm{-}\mathrm{var}})^c,
\nn\\
d_{p\textrm{-}\mathrm{var}} (\mathbf{y}, \hat{\mathbf{y}})
&\le
c (1 + \| f\|_{C_{\mathrm{b}}^{[p] +1}})^c
 (1+ \vertiii{\mathbf{x}}_{p\textrm{-}\mathrm{var}} + \vertiii{\hat{\mathbf{x}}}_{p\textrm{-}\mathrm{var}})^c
d_{p\textrm{-}\mathrm{var}} (\mathbf{x}, \hat{\mathbf{x}})
\nn
\end{align}
for all $\mathbf{x}, \hat{\mathbf{x}} \in G\Omega_p (\R^d)$.
\\
{\rm (ii)}~Suppose that there exists a constant $\kappa >0$ such that
\[
\max_{0\le j \le [p]+1} 
 \sup\{ |\nabla^j f (y)| \,:\, y\in\R^e, \, |y| \le R \}
= O (R^\kappa) \qquad \mbox{as $R \to\infty$,}
\]
where $O$ stands for the big Landau symbol.
Then, there exists
a constant $c' >0$ (independent of $\mathbf{x}$ and $\hat{\mathbf{x}}$) such that
\begin{align} 
\vertiii{\mathbf{y}}_{p\textrm{-}\mathrm{var}} 
&\le  
c' (1 + \|x\|_\infty^\kappa)^c (1+ \vertiii{\mathbf{x}}_{p\textrm{-}\mathrm{var}})^c,
\nn\\
d_{p\textrm{-}\mathrm{var}} (\mathbf{y}, \hat{\mathbf{y}})
&\le
c' (1 + \|x\|_\infty^\kappa + \|\hat{x}\|_\infty^\kappa)^c
 (1+ \vertiii{\mathbf{x}}_{p\textrm{-}\mathrm{var}} + \vertiii{\hat{\mathbf{x}}}_{p\textrm{-}\mathrm{var}})^c
d_{p\textrm{-}\mathrm{var}} (\mathbf{x}, \hat{\mathbf{x}})
\nn
\end{align}
for all $\mathbf{x}, \hat{\mathbf{x}} \in G\Omega_p (\R^d)$.
Here, we wrote $x_t:=\mathbf{x}^1_{0,t}$ for brevity
and $c$ is the same constant as in {\rm (ii)} above.
\end{proposition}

\begin{proof} 
{\rm (i)} is well-known (see \cite{lq, lcl} for instance). 
By a standard cut-off technique, we can easily obtain {\rm (ii)} from {\rm (i)}.
\end{proof}

\begin{remark} \label{rem.RPint1}
Instead of $f \colon \R^e \to \R^{e \times d}$ itself, 
we will often use 
$\mathrm{Id}_d \oplus f  \colon \R^e \to \R^{(d+e) \times d}$,
where $\mathrm{Id}_d$ stands for the identity map of $\R^d$.
Then, $\mathbf{x} \mapsto \int \tilde{f}(\mathbf{x}) d\mathbf{x}$
is the unique continuous extension of 
$x \mapsto (x, \, \int_0^{\cdot} f(x_s) dx_s)$.
\end{remark}

\begin{remark} \label{rem.RPint2}
The following operations can be viewed as special cases of 
RP integration. 
(In this remark, $m, n, l \in \N$.)
\\
{\rm (i)}~(Linear combination) Let $A, \hat{A}\in \R^{n\times n}$. Then, the map
 \begin{align*}
 G\Omega_p (\R^m\oplus \R^n \oplus \R^n)  \ni
 (\mathbf{z}, \mathbf{v}, \hat{\mathbf{v}})  &\mapsto
  (\mathbf{z}, \mathbf{v}, \hat{\mathbf{v}}, A\mathbf{v}+ \hat{A}\hat{\mathbf{v}})
  \in G\Omega_p (\R^m\oplus \R^n \oplus \R^n \oplus \R^n)
      \nn\\
      &\mapsto
  (\mathbf{z},  A\mathbf{v}+ \hat{A}\hat{\mathbf{v}})
  \in G\Omega_p (\R^m\oplus  \R^n)     
   \end{align*}
          is well-defined and locally Lipschitz continuous. 
            (The second arrow is the natural projection.)
            \\
            {\rm (ii)}~(Multiplication)  The map
             \begin{align*}
 G\Omega_p (\R^m\oplus \R^{l \times n} \oplus \R^n)  \ni
 (\mathbf{z}, \mathbf{J}, \mathbf{v})  &\mapsto
  (\mathbf{z}, \mathbf{J}, \mathbf{v},\,\mathbf{J}\mathbf{v})
  \in G\Omega_p (\R^m\oplus \R^{l \times n} \oplus \R^n \oplus \R^l)
   \end{align*}
            is well-defined and locally Lipschitz continuous. 
             (This is a unique continuous extension of the map  
             $(z_t, J_t, v_t)_{t \in [0,1]} \mapsto (z_t, J_t, v_t, J_t v_t)_{t \in [0,1]}$
              in the $1$-variational setting.
              Notice that $J_t v_t = \int_0^t J_s dv_s + \int_0^t (dJ_s )v_s$.)
                                      \end{remark}

\section{Estimates for RDE}\label{sec.Jac}

In this section we provide some useful estimates for RDEs.
In this paper RDEs are formulated in Lyons' original way 
(see \cite{lq, lcl}).
Throughout this section  $d, e\in\N$, $2 \le p <\infty$.
(Those are basically fixed unless otherwise stated
and we do not keep track of them.)
All results in Subsections 4.1 and 4.2 are basically known.
In Subsections 4.3 we estimate the difference of 
two solutions of a Jacobian RDE, which is the main result of this section.

\subsection{Standard RDE}
In this subsection we quickly review the standard
RDE theory of Lyons' type (see \cite{lq, lcl}).
The coefficient of RDE is $\sigma\in 
C_{\mathrm{b}}^{[p]+1} (\R^e, \R^{e \times d})$.
We write $K$ for the $C_{\mathrm{b}}^{[p]+1}$-norm of $\sigma$ for brevity.

For a driving RP $\mathbf{x} \in G\Omega_p (\R^d)$
and an initial value $a\in \R^e$,
we consider the following RDE in Lyons' sense:
\begin{equation}
dy_t = \sigma (y_t) dx_t, 
\qquad
\quad  y_0 =a.
\label{def.RDE_lyons}
\end{equation} 
(We do not use boldface letters when in the description of an RDE.
Instead, we use corresponding Italic letters.
Only in this section, we do not assume $a=0$.)
Its precise meaning is as follows: 
First, set $\sigma_a (y) := \sigma (a +y)$ and 
 \[
\tilde{\sigma}_a (z) \la z' \ra 
=
\begin{pmatrix}
{\rm Id} & 0 
\\ 
\sigma_a ( y ) &0
\end{pmatrix}
\begin{pmatrix}
x' \\ y'
\end{pmatrix},
\qquad
z=\begin{pmatrix}
x \\ y
\end{pmatrix},
 z'=\begin{pmatrix}
x' \\ y'
\end{pmatrix} \in \R^{d+e}= \R^d \oplus \R^e.
\]
Then, $\tilde{\sigma}_a \in 
C_{\mathrm{b}}^{[p]+1} (\R^{e+d}, \R^{(e+d) \times (e+ d)})$
and its $C_{\mathrm{b}}^{[p]+1}$ is bounded in $a$.
(We do not clearly distinguish row and column vectors in this paper.)
A solution of \eqref{def.RDE_lyons} is defined to be 
 $\mathbf{z}\in G\Omega_p (\R^{d+e})$
such that
 \begin{equation}\label{of2.def.rgh}
\mathbf{z}= \int \tilde{\sigma}_a (\mathbf{z}) d\mathbf{z}
\qquad
\mbox{with} \quad
\pi_1 \mathbf{z} =\mathbf{x},
\end{equation} 
where the integral is an RP integral of Lyons' sense 
and $\pi_1 \colon G\Omega_p (\R^{d+e}) \to G\Omega_p (\R^{d})$
is a natural projection onto ``the first component."
We sometimes call $\mathbf{y} := \pi_2 \mathbf{z} \in G\Omega_p (\R^{e})$
a solution, too. 
Here, $\pi_2 \colon G\Omega_p (\R^{d+e}) \to G\Omega_p (\R^{e})$
is a natural projection onto ``the second component."
(By replacing $[0,1]$ by a subinterval $[0,\tau]$ with 
 $0<\tau \le 1$, we can define a local solution of the RDE in the usual way.)
We often write $y_t :=y_0 + \mathbf{y}^1_{0,t}$.

\begin{remark} \label{rem.pairRP}
Sometimes it is convenient to write $(\mathbf{x}, \mathbf{y})$ 
 for $\mathbf{z}$. 
  It helps heuristic understanding.
   A drawback of this notation is that $(\mathbf{x}, \mathbf{y})$
     looks like a pair of two geometric RPs. 
       When we use this notation, we must keep in mind that 
         $(\mathbf{x}, \mathbf{y})$ belongs to $G\Omega_p (\R^{d+e})$
           (not to $G\Omega_p (\R^{d})\times G\Omega_p (\R^{e})$).
   \end{remark}

\begin{remark} \label{rem.vect_rde}
We often write $\sigma = [V_1, \ldots, V_d]$, where 
we view $V_j \in C_{\mathrm{b}}^{[p]+1} (\R^e, \R^e)$ for every $1\le j\le d$.
Then, RDE \eqref{def.RDE_lyons} can be written 
in an equivalently way as follows:
\[
dy_t = \sum_{j=1}^d V_j (y_t) dx_t^j, 
\qquad
\quad  y_0 =a.
\]
Here, the superindex ``$j$" of $dx^j$ stands for the coordinate 
 (not the level of an RP). 
When we use this notation, we sometimes view $V_i$'s as vector fields 
(i.e. first-order differential operators) on $\R^e$.
\end{remark}

Let us first recall Lyons' continuity theorem, 
which is the most important theorem in RP theory.
\begin{proposition} \label{RDE_welldef}
Let the situation be as above. 
  Then, for every $\mathbf{x}\in G\Omega_p (\R^{d})$ and $a\in\R^e$, a unique (global) solution 
   $\mathbf{z} \in G\Omega_p (\R^{d+e})$ exists. 
     Moreover, we have the following two estimates:
       \\
       {\rm (i)}~Suppose that,
        for a control function $\omega$ and $\tau\in (0,1]$, we have
        \[
         |\mathbf{x}_{s,t}^k | \le \omega(s,t)^{k/p}, 
          \qquad 
            0\le s \le t \le  \tau, \,\, 1 \le k \le [p].
        \] 
       Then, there exists a constant $C_1 >0$ 
       (which depends only on $\omega(0,\tau)$ and $K$) such that,
       \[
       |\mathbf{z}_{s,t}^k | \le C_1\omega(s,t)^{k/p}, 
          \qquad 
            0\le s \le t \le  \tau, \,\, 1 \le k \le [p].
               \]
                                   \\
                  {\rm (ii)}~Let $\mathbf{z}$ and $\hat{\mathbf{z}}$
 be a unique solution associated with  $(\mathbf{x}, a)$ and 
   $(\hat{\mathbf{x}}, \hat{a})$, respectively.  
    Suppose that,
        for a control function $\omega$, $\tau\in (0,1]$ and $\ve \ge 0$, we have
        \begin{eqnarray*}
         |\mathbf{x}_{s,t}^k | \vee  |\hat{\mathbf{x}}_{s,t}^k | \le \omega(s,t)^{k/p}, 
          \quad 
           |\mathbf{x}_{s,t}^k  -\hat{\mathbf{x}}_{s,t}^k |   \le \ve\omega(s,t)^{k/p}, 
            \quad  
               0\le s \le t \le  \tau, \,\, 1 \le k \le [p].
        \end{eqnarray*}
Then, there exists a constant $C_2 >0$ 
       (which depends only on $\omega(0,\tau)$ and $K$) such that,
         \[
               |\mathbf{z}_{s,t}^k -  \hat{\mathbf{z}}_{s,t}^k| 
                 \le C_2 (\ve + |a -\hat{a}|) \omega(s,t)^{k/p}, 
          \qquad 
            0\le s \le t \le  \tau, \,\, 1 \le k \le [p].
                     \]
                                   \end{proposition}

\begin{proof} 
This result is well-known. See \cite{lq, lcl}.
\end{proof}

For $x \in  \cC^{1\textrm{-}\mathrm{var}}_0 (\R^{d})$,
denote by $y \in  \cC^{1\textrm{-}\mathrm{var}}_a (\R^{d})$ the solution of Eq. \eqref{def.RDE_lyons} understood in the Riemann-Stieltjes (or Young) sense.
  Then, the solution of RDE \eqref{def.RDE_lyons}
     driven by $\mathbf{x} =S_p (x)$ coincides $S_p ((x, y_{\cdot} -a))$ and,
       in particular, $y_t = a + \mathbf{y}^{1}_{0,t}$. 
         (Recall that in Lyons' original formulation of RDEs, 
           the initial value of a solution must be adjusted.)
From this fact, we can easily see that 
a solution of the RDE is consistent in $p$.
Namely, if $\mathbf{z}$ solves RDE \eqref{def.RDE_lyons}
     driven by $\mathbf{x}$ in the $p$-variation setting, 
       then $\mathbf{Ext}_{p,r} (\mathbf{z})$ solves the same RDE
         driven by $\mathbf{Ext}_{p,r} (\mathbf{x})$ in the $r$-variation setting ($p \le r$).

Now we give a simple lemma, which enables us 
to obtain a global estimate of an RP from local ones.
\begin{lemma} \label{lem.loc_to_glo}
Let $n \in \N$,  $\ve \ge 0$,
$\mathbf{w}, \hat{\mathbf{w}} \in \Omega_p (\R^n)$
and $\omega\colon \triangle \to [0,\infty)$ be a control function.
Denote by 
$\cP =\{0 =\tau_0 <\tau_1 <\cdots <\tau_N= 1\}$ 
 a partition of $[0,1]$.
\\
{\rm (i)}~Suppose that 
\[
|\mathbf{w}^k_{s,t}| \le \omega (s,t)^{k/p}, \qquad 
\tau_{j-1}\le s \le t \le  \tau_j, \,\, 1 \le k \le [p]
\]
holds for all $1\le j\le N$.  Then, 
\begin{equation}\label{eq.0117-1}
|\mathbf{w}^k_{s,t}| \le N^{k(p-1)/p} \omega (s,t)^{k/p}, 
\qquad 
(s,t)\in \triangle, \,\, 1 \le k \le [p].
\end{equation}
{\rm (ii)}~
Suppose that both $\mathbf{w}$ and $\hat{\mathbf{w}}$
satisfy the assumption of {\rm (i)} above and, moreover,
\[
|\mathbf{w}^k_{s,t} - \hat{\mathbf{w}}^k_{s,t}| \le \ve \omega (s,t)^{k/p}, \qquad 
 \tau_{j-1}\le s \le t \le  \tau_j, \,\, 1 \le k \le [p]
\]
holds for all $1\le j\le N$.  Then, 
\begin{equation}\label{eq.0117-2}
|\mathbf{w}^k_{s,t} - \hat{\mathbf{w}}^k_{s,t}|  \le \ve N^{k(p-1)/p} \omega (s,t)^{k/p}, \qquad (s,t)\in \triangle, \,\, 1 \le k \le [p].
\end{equation}
 \end{lemma}

\begin{proof} 
Define $s=t_0 < t_1 <\cdots < t_M =t$ ($1 \le M \le N$)
  be such that $t_1, \ldots, t_{M-1}$ be  
   all $\tau_j$'s belonging to $(s,t)$.
    We prove \eqref{eq.0117-1} and  \eqref{eq.0117-2} 
      (with $N$ being replaced by $M$)
            by mathematical induction with respect to $M$.
              When $M=1$, there is nothing to prove since 
                $s$ and $t$ belong to the same subinterval in this case.
                 
Now we assume the case $M-1$ and will prove the case $M$.   
Chen's identity reads
       \begin{equation} \label{eq.0117-3}
       \mathbf{w}^k_{s,t} = 
          \mathbf{w}^k_{s,u} + \sum_{i=1}^{k-1} 
\mathbf{w}^{k-i}_{s,u}\otimes \mathbf{w}^{i}_{u,t}  + \mathbf{w}^k_{u,t},
\qquad 
1 \le k \le [p], \,\,
 s \le u \le t.
 \end{equation}
 We use \eqref{eq.0117-3} with $(s,u,t)= (t_0, t_{M-1}, t_M)$.  
Using the assumption of the induction, we see that
 \begin{align} 
|\mathbf{w}^k_{t_0,t_M} |
 &\le  (M-1)^{k(p-1)/p}\omega (t_0, t_{M-1})^{k/p} 
   \nn\\
   &\quad 
   + \sum_{i=1}^{k-1} (M-1)^{i(p-1)/p}\omega (t_0, t_{M-1})^{i/p} 
       \omega (t_{M-1}, t_M)^{(k-i)/p}
          + \omega (t_{M-1}, t_M)^{k/p}
          \nn\\ 
            &\le \bigl\{  
               (M-1)^{(p-1)/p}\omega (t_0, t_{M-1})^{1/p}   
               + 1^{(p-1)/p}\omega (t_{M-1}, t_M)^{1/p}   
                  \bigr\}^k
                      \nn\\ 
            &\le
               \{  M^{(p-1)/p}  \omega (t_0, t_{M})^{1/p}  \}^k,              
               \end{align}
                which is the desired inequality.  
                In the last inequality, we have used H\"older's inequality
                 and the superadditivity for $\omega$.
                Thus, we have shown \eqref{eq.0117-1}.
                
To prove \eqref{eq.0117-2}, we estimate the difference
of \eqref{eq.0117-3} for $\mathbf{w}^k$ and that for $\hat{\mathbf{w}}^k$.
Then, in the same way as above, 
 \begin{align} 
|\mathbf{w}^k_{t_0,t_M} - \hat{\mathbf{w}}^k_{t_0,t_M}| 
 &\le \ve (M-1)^{k(p-1)/p}\omega (t_0, t_{M-1})^{k/p} 
   \nn\\
   &\qquad 
   + 2\ve \sum_{i=1}^{k-1} (M-1)^{i(p-1)/p}\omega (t_0, t_{M-1})^{i/p} 
       \omega (t_{M-1}, t_M)^{(k-i)/p}
          \nn\\
   &\qquad           +\ve \omega (t_{M-1}, t_M)^{k/p}
          \nn\\ 
            &\le \ve \bigl\{  
               (M-1)^{(p-1)/p}\omega (t_0, t_{M-1})^{1/p}   
               + 1^{(p-1)/p}\omega (t_{M-1}, t_M)^{1/p}   
                  \bigr\}^k
                      \nn\\ 
            &\le
               \ve \{  M^{(p-1)/p}  \omega (t_0, t_{M})^{1/p}  \}^k.             
               \end{align}
                Thus, we have shown \eqref{eq.0117-2}, too.
                                (Note that the case $k=1$ is also covered in this argument. When $k=1$, the terms of the form $\sum_{i=1}^{k-1} (\cdots)$
                                 should be understood as $0$.) 
                                   This completes the proof of the lemma.
                     \end{proof}
 

Now we provide two corollaries of Proposition \ref{RDE_welldef}.
Recall that $N^p_{\beta} (\mathbf{x})$ was defined in \eqref{def_N}.
Recall that $N^p_{\beta} (\mathbf{x}) \le \beta^{-1} \vertiii{\mathbf{x}}_{p\textrm{-}\mathrm{var}}^p$.

\begin{corollary} \label{cor.0117-0}
Let $\mathbf{z}$ be a unique solution of RDE \eqref{def.RDE_lyons}
   associated with  $(\mathbf{x}, a)$. 
     Then, there exists  a constant $c=c(K)>0$ which depends 
       only on $K := \|\sigma\|_{C_{\mathrm{b}}^{[p]+1} }$ such that
   \begin{equation}\label{eq.0117-7}    
         \vertiii{\mathbf{z}}_{p\textrm{-}\mathrm{var}} 
         \le  c ( \vertiii{\mathbf{x}}_{p\textrm{-}\mathrm{var}} +\vertiii{\mathbf{x}}_{p\textrm{-}\mathrm{var}} ^p).
       \end{equation}
       More precisely, there exists a constant $c' >0$ (which is independent of 
        $\sigma, \mathbf{x}, a$) such that
          \begin{equation}\label{eq.0117-8}
           \vertiii{\Gamma_{A(K)} \mathbf{z}}_{p\textrm{-}\mathrm{var}} 
         \le  c' ( K\vertiii{\mathbf{x}}_{p\textrm{-}\mathrm{var}} +K^p \vertiii{\mathbf{x}}_{p\textrm{-}\mathrm{var}} ^p).
                   \end{equation}
                   Here, $A(K) \in \R^{(d+e) \times (d+e)}$ is a block matrix of the form 
 \[
 A(K) 
= 
\begin{pmatrix}
K \mathrm{Id}_d & 0 \\
 0& \mathrm{Id}_e
\end{pmatrix}
\]
and $\Gamma_{A(K)}$ is the generalized dilation.
In particular, there exists a constant $c'' >0$ (which is independnent of 
        $\sigma, \mathbf{x}, a$) such that 
          \begin{equation}\label{eq.0117-9}
           \vertiii{ \mathbf{z}}_{p\textrm{-}\mathrm{var}} 
         \le  c'' (1+K)^p( \vertiii{\mathbf{x}}_{p\textrm{-}\mathrm{var}} + \vertiii{\mathbf{x}}_{p\textrm{-}\mathrm{var}} ^p).
                   \end{equation}
                                      \end{corollary}

\begin{proof} 
Let $\{\tau_m : 1\le m \le N^p_{1} (\mathbf{x}) +1\}$ be the sequence that appears in the definition of $N^p_{1} (\mathbf{x})$.
Set $\omega (s,t) = \vertiii{\mathbf{x}}_{p\textrm{-}\mathrm{var}, [s,t]}^p$.
Then, $\omega ( \tau_{m-1}, \tau_m) \le 1$ for all $m$.
Proposition \ref{RDE_welldef} {\rm (i)} implies that there exists 
a constant $c_1>0$ such that
\[
 |\mathbf{z}_{s,t}^k | \le \{c_1\omega(s,t)\}^{k/p}, 
          \qquad 
             \tau_{m-1} \le s \le t \le  \tau_m, \,\, 1 \le k \le [p].
            \]
            on each subinterval.  
              Note that $c_1$ depends only on  $K$ (but not on $m$).  
            By  Lemma \ref{lem.loc_to_glo} {\rm (i)}, we have        
            \begin{align*}
            |\mathbf{z}^k_{s,t}| &\le ( N^p_{1} (\mathbf{x}) +1)^{k(p-1)/p} 
              \{c_1    \omega (s,t)\}^{k/p}
            \\
            &\le 
            \{ c_1   
              ( N^p_{1} (\mathbf{x}) +1)^{p-1} \omega (s,t) \}^{k/p}, 
\quad 
(s,t)\in \triangle, \,\, 1 \le k \le [p].
            \end{align*}      
            One can easily see from this that
            \[
            \|\mathbf{z}^k\|_{p/k\textrm{-}\mathrm{var}} 
            \le 
                  \{ c_1  (\vertiii{\mathbf{x}}_{p\textrm{-}\mathrm{var}}^p +1)^{(p-1)/p}  \omega (0,1)^{1/p} \}^{k}
\le
\{ c_2  (\vertiii{\mathbf{x}}_{p\textrm{-}\mathrm{var}}^p +  \vertiii{\mathbf{x}}_{p\textrm{-}\mathrm{var}})   \}^{k},
\]
where $c_2>0$ is a certain constant depending only on $K$.
Thus, we obtained \eqref{eq.0117-7}.     

By using the standard scaling argument,    
we can show \eqref{eq.0117-8} from \eqref{eq.0117-7} with $K=1$
as follows.
If $\mathbf{z} =(\mathbf{x}, \mathbf{y})$ solves
RDE \eqref{def.RDE_lyons} associated with $\sigma$ and $(\mathbf{x}, a)$, 
then 
$\Gamma_{A(K)} \mathbf{z} \,\,``=(K\mathbf{x}, \mathbf{y})"$
associated with $\sigma/K$ and $(K\mathbf{x}, a)$.
When $\mathbf{x}$ is a natural lift of  a path of finite $1$-variation,  
 this fact can be checked easily.
 When $\mathbf{x}$ is a general geometric RP, 
 approximate it by a sequence of such nice paths. 
 Since the norm of $\sigma/K$ equals $1$ (unless $\sigma$ is trivial), 
 we can use \eqref{eq.0117-7} with $K=1$ to prove 
 \eqref{eq.0117-8}  with $c' = c(1)$.
 
 Finally, \eqref{eq.0117-9} follows immediately from \eqref{eq.0117-8}.
     \end{proof}


%
\begin{corollary} \label{cor.0117-1}
Let $\mathbf{z}$ and $\hat{\mathbf{z}}$
 be a unique solution of RDE \eqref{def.RDE_lyons}
   associated with  $(\mathbf{x}, a)$ and 
   $(\hat{\mathbf{x}}, \hat{a})$, respectively.  
      Then, for every $\beta >0$, there exists $c>0$ such that
    \begin{align*} 
         d_{p\textrm{-}\mathrm{var}} (\mathbf{z}, \hat{\mathbf{z}})
                  &\le 
          c (\vertiii{\mathbf{x}}_{p\textrm{-}\mathrm{var}}^p + 
 \vertiii{\hat{\mathbf{x}}}_{p\textrm{-}\mathrm{var}}^p +1)
          \\
          & \qquad \times
          \exp \bigl[
            c (N^p_{\beta} (\mathbf{x})+ N^p_{\beta} (\hat{\mathbf{x}}))
            \bigr]
          (d_{p\textrm{-}\mathrm{var}} (\mathbf{x}, \hat{\mathbf{x}}) +|a-\hat{a}|).
         \end{align*}
         Here, the constant $c$ depends only on $\beta$ and $K$.
   \end{corollary}

\begin{proof} 
Let $\{\tau_m\}$ and $\{\hat{\tau}_m\}$ be the sequence that appears
in the definition of $N^p_{\beta} (\mathbf{x})$ and 
 $N^p_{\beta} (\hat{\mathbf{x}})$, respectively.
 Let $0=u_0 <u_1 <\cdots < u_M =1$ be 
 all elements of $\{\tau_m\}\cup \{\hat{\tau}_m\}$ in increasing order.
 Clearly, $M \le N^p_{\beta} (\mathbf{x})+N^p_{\beta} (\hat{\mathbf{x}}) +2$.
  We write $y_t = a + \mathbf{y}^1_{0,t}$ and  
 $\hat{y}_t= \hat{a} + \hat{\mathbf{y}}^1_{0,t}$.  
 In this proof, $c_i~(i=1,2,\dots)$ are positive constants which depend
 only on $\beta$ and $K$.

Define a new control function $\omega$ by 
\begin{equation}\label{def.diff_omega}
\omega (s,t) := \vertiii{\mathbf{x}}_{p\textrm{-}\mathrm{var}, [s,t]}^p + 
 \vertiii{\hat{\mathbf{x}}}_{p\textrm{-}\mathrm{var}, [s,t]}^p 
   + 
   \sum_{1 \le i \le [p]} \ve^{-p/i}
      \|  \mathbf{x}^i - \hat{\mathbf{x}}\|_{p/i\textrm{-}\mathrm{var}, [s,t]}^{p/i},
        \end{equation}
where we set $\ve := d_{p\textrm{-}\mathrm{var}} (\mathbf{x}, \hat{\mathbf{x}})$
for simplicity.
By construction, $\omega (u_{j-1}, u_j) \le 2\beta +[p]$  for each $j$.
Moreover, $\omega$ satisfies
   the assumptions of Proposition \ref{RDE_welldef}  {\rm (ii)}
     on each subinterval $[u_{j-1}, u_j]$. 
       Hence, there exist $c_1$ and $c_2$ such that, for all $1 \le j\le M$,
         $u_{j-1}\le s \le t \le u_{j}$ and $1 \le k \le [p]$, we have
                \begin{align}
              |\mathbf{z}_{s,t}^k |\vee |\hat{\mathbf{z}}_{s,t}^k| 
                 &\le \{ c_1 \omega(s,t)\}^{k/p}, 
                   \quad
                     |\mathbf{z}_{s,t}^k -  \hat{\mathbf{z}}_{s,t}^k| 
                 \le c_2 (\ve + |y_{u_{j-1}} -\hat{y}_{u_{j-1}}|) \omega(s,t)^{k/p}.
          \label{eq.0117-5}
       \end{align}   

Write $D:= c_2 (2\beta +[p])^{1/p}$ for simplicity.
From the second inequality in \eqref{eq.0117-5},  we see that      
         \begin{align}
         |y_{u_j} - \hat{y}_{u_j}| 
           &\le  
             |y_{u_j} - \hat{y}_{u_j}| +   
               |\mathbf{y}_{u_{j-1}, u_j}^1 -  \hat{\mathbf{y}}_{u_{j-1}, u_j}^1| 
                \nn\\
                  &\le 
                    D \ve +(D+1)|y_{u_{j-1}} - \hat{y}_{u_{j-1}}|
              \nn\\
                 &\le 
              D \ve +(D+1) \{
               D \ve +(D+1)|y_{u_{j-2}} - \hat{y}_{u_{j-2}}| \}
                  \nn\\
                      &\le 
                         D\{ 1+ (D+1)\}\ve + (D+1)^2|y_{u_{j-2}} - \hat{y}_{u_{j-2}}|
                         \nn\\         
                             &\le 
                              D\{ 1+ (D+1)+\cdots +(D+1)^{j-1}\}\ve +(D+1)^j 
                               |y_{u_{0}} - \hat{y}_{u_{0}}|
                              \nn\\
                              &\le 
                                 \{(D+1)^j -1 \} \ve  +(D+1)^j  |a-\hat{a}|.                             
                                                                                           \label{ineq.0123-2}
                                                                                              \end{align} 
Note that we have essentially seen that
\[
\|y -\hat{y}\|_\infty \le (D+1)^{M}  (\ve +   |a-\hat{a}|). 
\]
Putting \eqref{ineq.0123-2} back into \eqref{eq.0117-5},  we  obtain that                 
 \[
     |\mathbf{z}_{s,t}^k -  \hat{\mathbf{z}}_{s,t}^k| 
                 \le c_2 (D+1)^M  \{\ve +|a-\hat{a}|\} \omega(s,t)^{k/p},
                    \quad u_{j-1}\le s \le t \le u_{j},  1 \le k \le [p]
                    \]                
                         for all $1 \le j\le M$.
                         
Combining this with the first inequality of \eqref{eq.0117-5}
and using Lemma \ref{lem.loc_to_glo} with the control $c_1\omega$,
  we see that
  \[
     |\mathbf{z}_{s,t}^k -  \hat{\mathbf{z}}_{s,t}^k| 
                 \le c_3 M^{k(p-1)/p}                 
                  e^{M\log (D+1)}  \{\ve +|a-\hat{a}|\} \omega(s,t)^{k/p},
                    \quad  (s,t) \in \triangle,  1 \le k \le [p]
  \]
    for some $c_3 >0$.       Finally,       
      by choosing a suitable constant $c > c_3 \vee \log (D+1)$,
        we can prove the desired estimate.  
        Noting that
        $\omega (0,1) \le \vertiii{\mathbf{x}}_{p\textrm{-}\mathrm{var}}^p + 
 \vertiii{\hat{\mathbf{x}}}_{p\textrm{-}\mathrm{var}}^p  + [p]$,  we finish the proof.
                          \end{proof}

\subsection{Jacobian RDE}

 In this subsection 
the coefficient is $\sigma\in 
C_{\mathrm{b}}^{[p]+2} (\R^e, \R^{e \times d})$.
We write $K'$ for the $C_{\mathrm{b}}^{[p]+2}$-norm of $\sigma$ for brevity.
(Since we study the derivative of the original RDE, 
one more differentiability is assumed.)
 We write $\mathcal{V}:=\R^e \oplus \R^{e \times e} 
\oplus \R^{e \times e}$ for simplicity.

For a driving RP $\mathbf{x} \in G\Omega_p (\R^d)$
and an initial value $(a, A, B) \in \mathcal{V}$,
we consider the following system of RDEs in Lyons' sense:
\begin{align} 
dy_t &= \sigma (y_t) dx_t, 
\qquad
\quad  y_0 =a,
\label{def.JRDE_1} 
\\
dJ_t &= \nabla\sigma (y_t)
  \langle  J_t, dx_t \rangle, 
\qquad
\quad  J_0 =A,
\label{def.JRDE_2} 
\\
dK_t &= - K_t \cdot  \nabla\sigma (y_t)
  \langle  \bullet, dx_t \rangle, 
\qquad
\quad  K_0 =B.
\label{def.JRDE_3} 
\end{align} 
Here, the right hand side  of \eqref{def.JRDE_3}
is (the minus of) the matrix multiplication 
of $K_t \in \R^{e \times e} $ and $ \nabla\sigma (y_t)\langle  \bullet, dx_t \rangle\in \R^{e \times e}$.
If we set 
\begin{equation} \label{def.M}
M_t := \int_0^t \nabla\sigma (y_s)\langle  \bullet, dx_s \rangle
\end{equation} 
as an $\R^{e \times e}$-valued integral,  
 \eqref{def.JRDE_2}--\eqref{def.JRDE_3} read:
  \begin{align} 
dJ_t &=  (dM_t)\cdot  J_t, 
\qquad
\quad  J_0 =A,
\label{def.JRDE_4} 
\\
dK_t &= - K_t  \cdot dM_t, 
\qquad
\quad  K_0 =B.
\label{def.JRDE_5} 
\end{align} 
A  solution of the system of RDEs
\eqref{def.JRDE_1}--\eqref{def.JRDE_3} 
is denoted by $(\mathbf{x}, \mathbf{y}, \mathbf{J}, \mathbf{K})$, 
which belongs to $G\Omega_p (\R^d \oplus \mathcal{V})$.
We call it a solution of {\rm (JacRDE)} driven by $\mathbf{x}$
with the initial value  $(a,A,B)$ for simplicity.
As before, we write 
\[
(y_t, J_t, K_t) :=(y_0 + \mathbf{y}^1_{0,t}, J_0 + \mathbf{J}^1_{0,t}, 
K_0 + \mathbf{K}^1_{0,t}).
\]
As is well-known,
when $J_0=\mathrm{Id}_e=K_0$,  we always have $J_t K_t \equiv \mathrm{Id}_e$.

\begin{remark} \label{rem.equivRDE}
Equivalently, \eqref{def.JRDE_1}--\eqref{def.JRDE_3}
can be written as in Remark \ref{rem.vect_rde}  in the following way:
\begin{align} 
dy_t &= \sum_{j=1}^d V_j (y_t) dx_t^j, 
\qquad
\quad  y_0 =a,
\nn
\\
dJ_t &= \sum_{j=1}^d \nabla V_j (y_t) J_t dx_t^j,
\qquad
\quad  J_0 =A,
\nn
\\
dK_t &= -  \sum_{j=1}^d K_t \nabla V_j (y_t) dx_t^j,
\qquad
\quad  K_0 =B.
\nn
\end{align} 
Here, $\nabla V_j$ is viewed as an $\R^{e \times e}$-valued function 
and the superindex ``$j$" of $dx^j$ stands for the coordinate 
 (not the level of an RP).
 \end{remark}

When viewed as a $\mathcal{V}$-valued RDE, the system of RDEs
\eqref{def.JRDE_1}--\eqref{def.JRDE_3} 
has a $C^{[p]+1}$-coefficient, 
which is not of $C_{\mathrm{b}}^{[p]+1}$ (due to the linearly 
growing property in the ``$(J, K)$-direction").
Hence, it always has a unique local solution.
But, it is not obvious whether a global solution exists and 
Lyons' continuity theorem holds or not.
    However, it is known that a unique global solution of {\rm (JacRDE)} exists 
    for every $\mathbf{x} \in G\Omega_p (\R^d)$ and $(a,A,B)$.
(See \cite[Section 10.7]{fvbk} for example.)
       Therefore, by a standard cut-off argument,  we can show that 
       Lyons' continuity theorem holds, namely,
        the map 
        \[ 
          (\mathbf{x} , (a, A, B)) \mapsto (\mathbf{x}, \mathbf{y}, \mathbf{J}, \mathbf{K})
                  \]
                   is locally Lipschitz continuous.
              In this subsection we recall estimates for 
                $(\mathbf{x}, \mathbf{y}, \mathbf{J}, \mathbf{K})$
                  which will be used in the next subsection.


As one can easily guess from \eqref{def.JRDE_4}--\eqref{def.JRDE_5},
$J$ and $K$ has right- and left-invariance, respectively.
Denote by $\tilde{\Gamma}_{A,B}$ the generalized dilation 
by the following linear map from $\R^d \oplus \mathcal{V}$ to itself:
\begin{equation}\label{def.Gamma_tilde}
\R^d \oplus \mathcal{V} \,\, \ni
(x,y, J, K)  \mapsto (x,y, JA, BK) \,\,  \in \R^d \oplus \mathcal{V}.
\end{equation}

\begin{lemma} \label{lem.invJK}
Let $\mathbf{x}\in G\Omega_p (\R^d)$ and $a \in \R^e$, 
$A, B \in \R^{e \times e}$.
Denote by $(\mathbf{x}, \mathbf{y}, \mathbf{J}, \mathbf{K})$
a unique 
solution of {\rm (JacRDE)} driven by $\mathbf{x}$
starting at $(a, \mathrm{Id}_e, \mathrm{Id}_e)$ on a certain time interval $[0, \tau]$. 
Then, 
$\tilde{\Gamma}_{A,B}(\mathbf{x}, \mathbf{y}, \mathbf{J}, \mathbf{K})$
is a unique 
solution of {\rm (JacRDE)} driven by $\mathbf{x}$
starting at $(a, A, B)$ on $[0, \tau]$. 
\end{lemma}

\begin{proof} 
Consider the case $\mathbf{x} = S_p (x)$ for some 
$x \in \cC^{1\textrm{-}\mathrm{var}}_0 (\R^d)$.
Then, $(y_t, J_t, K_t)$ coincides with the unique solution of 
  \eqref{def.JRDE_1}--\eqref{def.JRDE_3} 
    starting at $(a, \mathrm{Id}_e, \mathrm{Id}_e)$
       understood as a system of Riemann-Stietjes (or Young) ODEs.
        In this case, we can easily see that $t \mapsto (y_t, J_t A, BK_t)$ 
          solves the same equation with the initial value 
            being replaced by $(a, A, B)$.  
              Since the natural lift of $t \mapsto (x_t, y_t, J_t A, BK_t)$ 
                coincides with $\tilde{\Gamma}_{A,B}(\mathbf{x}, \mathbf{y}, \mathbf{J}, \mathbf{K})$, our assertion is true in this case.

For a general geometric RP $\mathbf{x}$,  
   we first approximate it by $\{S_p (x_n)\}_{n \in \N}$, for a
    certain $\{x_n\}_{n\in\N} \subset \cC^{1\textrm{-}\mathrm{var}}_0 (\R^d)$ 
      and then use Lyons' continuity theorem.
       (Note that even though the coefficient of this RDE is not $C_{\mathrm{b}}^{[p]+1}$,
        the  continuity theorem still works for local solutions 
         thanks to a standard cut-off technique.)
 \end{proof}


\begin{lemma} \label{lem.getout}
Let $\mathbf{x}\in G\Omega_p (\R^d)$, $a \in \R^e$ and $\tau \in (0, 1]$.
Denote by $(\mathbf{x}, \mathbf{y}, \mathbf{J}, \mathbf{K})$
a unique local solution defined on $[0,\tau]$
of {\rm (JacRDE)} driven by $\mathbf{x}$
starting at $(a, \mathrm{Id}_e, \mathrm{Id}_e)$. 
Then, there exists a positive constant $\beta_0$ 
depending only on $K'$ with the following property:  

If $\vertiii{\mathbf{x}}_{p\textrm{-}\mathrm{var}, [0,\tau]}^p \le \beta_0$
         holds, then $\|J\|_{\infty, [0, \tau]} \vee \|K\|_{\infty, [0, \tau]} 
          \le 3 |\mathrm{Id}_e|$.
         Here, $\|\cdot\|_{\infty, [0, \tau]}$ stands for the usual sup-norm over
            the interval $[0, \tau]$.
\end{lemma}

\begin{proof} 
Let $\chi\colon \R^{e\times e} \to [0,1]$ be a smooth, 
function with compact support such that $\chi (J) =1$ 
  if $|J| \le 3 |\mathrm{Id}_e|$.
   We cut-off the coefficient of {\rm (JacRDE)} by using this $\chi$.
     (Just replace $J_t$ and $K_t$
      on the right hand side of \eqref{def.JRDE_2}  and \eqref{def.JRDE_3}
            by $\chi (J_t)J_t$ and $\chi (K_t)K_t$, respectively.)
             Then, the new system of 
                RDEs has a $C_{\mathrm{b}}^{[p]+1}$-coefficient.
               So, we can use 
               Proposition \ref{RDE_welldef} {\rm (i)}. 
                    It implies that if $\omega(0,\tau)$ is small enough, 
                      the $J$-component and $K$-component
                        of the first level of the new system 
                        do not get out of the ball of radius $3 |\mathrm{Id}_e|$. 
                         (Hence, this solution of the truncated system  
                           also solves the original one.)
             \end{proof}


\begin{lemma} \label{lem.glob_est}
Let $\mathbf{x}\in G\Omega_p (\R^d)$ and $a \in \R^e$.
Then, there exists a unique global solution  
$(\mathbf{x}, \mathbf{y}, \mathbf{J}, \mathbf{K})$ 
of {\rm (JacRDE)} driven by $\mathbf{x}$
starting at $(a, \mathrm{Id}_e, \mathrm{Id}_e)$. 
Moreover, for every $\beta \in (0, \beta_0]$, there exists a 
constant $C>0$ depending only on $\beta$ and $K'$ such that
\begin{align}
\|J\|_{\infty} \vee \|K\|_{\infty} 
&\le
C \exp \bigl[ C N^p_{\beta} (\mathbf{x}) \bigr],
 \label{ineq.0119-1}\\
 \vertiii{(\mathbf{x}, \mathbf{y}, \mathbf{J}, \mathbf{K})
      }_{p\textrm{-}\mathrm{var}}
          &\le C \exp \bigl[ C N^p_{\beta} (\mathbf{x}) \bigr]
            ( \vertiii{\mathbf{x}}_{p\textrm{-}\mathrm{var}} + \vertiii{\mathbf{x}}_{p\textrm{-}\mathrm{var}} ^p).
                         \label{ineq.0119-2}
                         \end{align}
                        Here, $\beta_0 >0$ is the constant that appears in 
                          Lemma \ref{lem.getout}.
\end{lemma}

\begin{proof} 
Let $\{\tau_m : 0\le m \le N^p_{\beta} (\mathbf{x})+1\}$ be the sequence that appears
in the definition of $N^p_{\beta} (\mathbf{x})$.

Denote by 
$(\mathbf{x} (m), \mathbf{y}(m), \mathbf{J}(m), \mathbf{K}(m))$ 
be a solution of {\rm (JacRDE)} driven by $\mathbf{x}$
on the subinterval $[\tau_{m-1}, \tau_{m}]$
with the initial condition
$(y_{\tau_{m-1}}, \mathrm{Id}_e, \mathrm{Id}_e)$ at time $\tau_{m-1}$.
This satisfies the estimate in Lemma \ref{lem.getout}.
We see from Lemma \ref{lem.invJK}  that
\[
\tilde{\Gamma}_{ J_{\tau_{m-1}}, K_{\tau_{m-1}}}
(\mathbf{x} (m), \mathbf{y}(m), \mathbf{J}(m), \mathbf{K}(m))
\]
coincides with $(\mathbf{x}, \mathbf{y}, \mathbf{J}, \mathbf{K})$ 
on $[\tau_{m-1}, \tau_{m}]$.
By concatenating them all, we obtain a global solution. 
In particular,  we have
\[
 \sup_{ \tau_{m-1}\le t\le \tau_{m}} |J_t| \le |J_{\tau_{m-1}}| 
 \cdot 3 |\mathrm{Id}_e| 
     \le   (3 |\mathrm{Id}_e|)^m
 \]
 for all $m$. $K$ satisfies the same estimate.
So,  \eqref{ineq.0119-1} is satisfied if we take
   $C =(3 |\mathrm{Id}_e|)\vee \log (3 |\mathrm{Id}_e|) =3 |\mathrm{Id}_e|$.

Next, we prove \eqref{ineq.0119-2} by using the standard cut-off technique.
For every $R \ge 1$, we can find a smooth, compactly-supported
function $\chi_R \colon \R^{e\times e} \to [0,1]$ with the following propeties:
{\rm (i)} $\chi_R \equiv 1$ on the ball of radius $R$ centered at $0$.
{\rm (ii)} $\chi_R$ vanishes outside the ball of radius $2R$ centered at $0$.
{\rm (iii)} The $C_{\mathrm{b}}^{[p]+1}$-norm of 
$\chi_R$ is bounded in $R$.

Using $\chi_R$ (instead of $\chi$), we cut-off the coefficient of 
  {\rm (JacRDE)} in the same way as in Lemma \ref{lem.getout}.  
    Then, the $C_{\mathrm{b}}^{[p]+1}$-norm of the new coefficient 
      is bounded by $C' R$ for some constant $C' >0$ independent of $R$.     
 For a given $\mathbf{x}$, we use this cut-off with $R= C \exp \bigl[ C N^p_{\beta} (\mathbf{x}) \bigr]$. 
 By applying \eqref{eq.0117-9} in Corollary \ref{cor.0117-0} 
  with $K=C'C \exp \bigl[ C N^p_{\beta} (\mathbf{x}) \bigr]$,
 we obtain \eqref{ineq.0119-2} (after adjusting the constant $C$).
 \end{proof}

\begin{remark}
It is quite important  that the positive constants 
that appear in the estimates in this and previous subsection 
do not depend on the initial value $a$.
\end{remark}
     
\subsection{Difference of two solutions of Jacobian RDE}
\label{subsec.43}

In this subsection, we continue to work in the setting of the previous 
subsection. The constant $\beta_0 >0$ is the one that appeared in Lemma \ref{lem.getout}. (It depends only on $K'$.)

\begin{lemma} \label{lem.diff_loc_Jac}     
Let $\mathbf{x}, \hat{\mathbf{x}}
\in G\Omega_p (\R^d)$ and $a, \hat{a} \in \R^e$, 
Denote by $(\mathbf{x}, \mathbf{y}, \mathbf{J}, \mathbf{K})$
a unique solution of {\rm (JacRDE)} driven by $\mathbf{x}$
starting at $(a, \mathrm{Id}_e, \mathrm{Id}_e)$.
Also, denote by $(\hat{\mathbf{x}}, \hat{\mathbf{y}}, 
  \hat{\mathbf{J}}, \hat{\mathbf{K}})$ one driven by $\hat{\mathbf{x}}$
starting at $(\hat{a}, \mathrm{Id}_e, \mathrm{Id}_e)$. 
Suppose further that,
        for a control function $\omega$, $\tau\in (0,1]$ and $\ve \ge 0$, we have $\vertiii{\mathbf{x}}_{p\textrm{-}\mathrm{var}, [0,\tau]}^p
\vee  \vertiii{\hat{\mathbf{x}}}_{p\textrm{-}\mathrm{var}, [0,\tau]}^p\le \beta_0$
                          and 
                                  \begin{eqnarray*}
         |\mathbf{x}_{s,t}^k | \vee  |\hat{\mathbf{x}}_{s,t}^k | \le \omega(s,t)^{k/p}, 
          \quad 
           |\mathbf{x}_{s,t}^k  -\hat{\mathbf{x}}_{s,t}^k |   \le \ve\omega(s,t)^{k/p}, 
            \quad  
               0\le s \le t \le  \tau, \,\, 1 \le k \le [p].
                \end{eqnarray*}
               Then, there exist constants $C_3, C_4 >0$
       (which depends only on $\omega (0,\tau)$, $K'$) such that
         \begin{align*}
         |(\mathbf{x}, \mathbf{y}, \mathbf{J}, \mathbf{K})_{s,t}^k |
                  \vee  |(\hat{\mathbf{x}}, \hat{\mathbf{y}}, 
  \hat{\mathbf{J}}, \hat{\mathbf{K}})_{s,t}^k|                        
   &\le C_3\omega(s,t)^{k/p},
                  \\
               |(\mathbf{x}, \mathbf{y}, \mathbf{J}, \mathbf{K})_{s,t}^k 
               -  (\hat{\mathbf{x}}, \hat{\mathbf{y}}, 
  \hat{\mathbf{J}}, \hat{\mathbf{K}})_{s,t}^k| 
               &\le C_4 (\ve + |a -\hat{a}|) \omega(s,t)^{k/p}
                     \end{align*}
hold for all $0\le s \le t \le  \tau$ and $1 \le k \le [p]$.
                          \end{lemma}

 \begin{proof}   
By the assumption  and Lemma \ref{lem.getout}, 
the $J$- and  $K$-component of the solutions  do not 
get out of the the centered ball of radius $3 |\mathrm{id}_e|$.
Then, we can easily show  
this lemma by Lemma \ref{lem.getout} and 
Proposition \ref{RDE_welldef}  {\rm (ii)}.
(We use the same cut-off as in the proof of  Lemma \ref{lem.getout}
and then apply Proposition \ref{RDE_welldef} for the new RDE.)
  \end{proof}

 
For the rest of this section, we use the following notation.
For $(\mathbf{x}, a) \in G\Omega_p (\R^d)\times  \R^e$,
we denote by $(\mathbf{x}, \mathbf{y}, \mathbf{J}, \mathbf{K})$
a unique global solution of {\rm (JacRDE)} driven by $\mathbf{x}$
starting at $(a, \mathrm{Id}_e, \mathrm{Id}_e)$.
Also, we denote by 
$(\hat{\mathbf{x}}, \hat{\mathbf{y}}, 
  \hat{\mathbf{J}}, \hat{\mathbf{K}})$
a unique global solution  associated with $(\hat{\mathbf{x}}, \hat{a})$.
 
\begin{lemma} \label{lem.diff_1st_Jac}     
Let the notation be  as above.
Then, for every $\beta \in (0,\beta_0]$, there exists $c>0$ such that
\[
\|J -\hat{J}\|_{\infty} + \|K-\hat{K}\|_{\infty}
 \le 
      c \exp \bigl[
            c (N^p_{\beta} (\mathbf{x})+ N^p_{\beta} (\hat{\mathbf{x}}))
            \bigr]
          (d_{p\textrm{-}\mathrm{var}} (\mathbf{x}, \hat{\mathbf{x}}) +|a-\hat{a}|)
          \]
          for all $(\mathbf{x}, a),  (\hat{\mathbf{x}}, \hat{a}) 
\in G\Omega_p (\R^d)\times  \R^e$.
                    Here, the constant $c$ depends only on $\beta$ and $K'$.          
\end{lemma}

 \begin{proof}   
 In this proof, $c_i~(i=1,2,\dots)$ are positive constants which depends 
 only on $\beta$ and $K'$.
 Let $\{\tau_m\}$ and $\{\hat{\tau}_m\}$ be the sequence that appears
in the definition of $N^p_{\beta} (\mathbf{x})$ and 
 $N^p_{\beta} (\hat{\mathbf{x}})$, respectively.
 Let $0=u_0 <u_1 <\cdots < u_M =1$ be 
 all elements of $\{\tau_m\}\cup \{\hat{\tau}_m\}$ in increasing order.
 Clearly, 
 \begin{equation}\label{eq.0121-2}
  M \le N^p_{\beta} (\mathbf{x})+N^p_{\beta} (\hat{\mathbf{x}}) +2.
  \end{equation}

Define $\omega$ by \eqref{def.diff_omega} in the proof 
of Corollary \ref{cor.0117-1}.
Then, $\omega (u_{j-1}, u_j) \le 2\beta +[p]$  for each $j$.
Moreover,  it was essentially proved there that 
$\delta^\prime := \,\,  \| y -\hat{y}\|_{\infty}$ satisfies 
\begin{equation}\label{eq.0121-3}
\delta^\prime  
\le c_1^{M+1} (\ve + |a- \hat{a}|),
\end{equation}
where $c_1 (:=D+1 >1)$ is  the constant and $\ve := d_{p\textrm{-}\mathrm{var}} (\mathbf{x}, \hat{\mathbf{x}})$
as before.

Denote by 
$(\mathbf{x} (j), \mathbf{y}(j), \mathbf{J}(j), \mathbf{K}(j))$ 
be a solution of {\rm (JacRDE)} driven by $\mathbf{x}$
on the subinterval $[\tau_{j-1}, \tau_{j}]$
with the initial value
$(y_{\tau_{j-1}}, \mathrm{Id}_e, \mathrm{Id}_e)$ at the time $\tau_{j-1}$.
Recall that
\[
\tilde{\Gamma}_{ J_{\tau_{j-1}}, K_{\tau_{j-1}}}
(\mathbf{x} (j), \mathbf{y}(j), \mathbf{J}(j), \mathbf{K}(j))
\]
coincides with $(\mathbf{x}, \mathbf{y}, \mathbf{J}, \mathbf{K})$ 
on $[\tau_{j-1}, \tau_{j}]$.
(Of course, $(\hat{\mathbf{x}} (j), \hat{\mathbf{y}} (j), 
  \hat{\mathbf{J}}(j), \hat{\mathbf{K}}(j))$ is defined in the same way 
  and also satisfies this property.) 
   We can apply to Lemma \ref{lem.diff_loc_Jac} to them on each subinterval:
   \begin{align*}
                  |\mathbf{J}(j)_{\tau_{j-1}, t} -  \hat{\mathbf{J}}(j)_{\tau_{j-1}, t}| 
                    +  |\mathbf{K}(j)_{\tau_{j-1}, t} -  \hat{\mathbf{K}}(j)_{\tau_{j-1}, t}|                   
               &\le C_3 (\ve + \delta^\prime) \omega(\tau_{j-1},t)^{1/p},
                  \,\, t\in [\tau_{j-1}, \tau_{j}]
                                       \end{align*}
                                       and, in particular,
\begin{align*}
 |J(j)_t -\hat{J}(j)_t|                       + |K(j)_t -\hat{K}(j)_t|                  
               &\le c_2 (\ve + \delta^\prime),
                  \qquad t\in [\tau_{j-1}, \tau_{j}]                                       \end{align*}
                     for all  $1 \le j \le M$.
                                       Here, we set $c_2 :=C_3(2\beta +[p])^{1/p}$ and used  the condition that 
                                       the initial values are the same.

We will prove by mathematical induction that 
\begin{equation}\label{eq.0121-1}
 |J_t -\hat{J}_t|                       \vee |K_t -\hat{K}_t|                  
               \le   j (3 |\mathrm{Id}_e|)^{j-1}  c_2 (\ve + \delta^\prime),
                  \qquad t\in [\tau_{j-1}, \tau_{j}].   
                  \end{equation}
for all $1 \le j \le M$.
(We only estimate $K-\hat{K}$ since we can estimate 
  $J-\hat{J}$ essentially in the same way.)
    We have already showed \eqref{eq.0121-1} for $j=1$.
      Suppose that \eqref{eq.0121-1} holds for $1, 2, \ldots, j-1$.    
           Note that 
           \[
           K_t = K_{\tau_{j-1}} K(j)_t = \{K(1)_{\tau_{1}}\cdots  K(j-1)_{\tau_{j-1}}\}K(j)_t,
           \qquad  t\in [\tau_{j-1}, \tau_{j}].                                            
           \]
By Lemma \ref{lem.getout}, we have $|K(j)_t| \le 3 |\mathrm{Id}_e|$ and
   $|K_{\tau_{j-1}}| \le (3 |\mathrm{Id}_e|)^{j-1}$. 
Then,  we can easily see that, for $t\in [\tau_{j-1}, \tau_{j}]$,
  \begin{align*}
     |K_t -\hat{K}_t| 
       &\le
         |K_{\tau_{j-1}} K(j)_t - \hat{K}_{\tau_{j-1}} \hat{K}(j)_t|
          \\
            & \le 
              |K_{\tau_{j-1}} -\hat{K}_{\tau_{j-1}}| \cdot |K(j)_t|
               + |\hat{K}_{\tau_{j-1}}| \cdot |K(j)_t - \hat{K}(j)_t|
                 \\
                 &\le 
                     (j-1) (3 |\mathrm{Id}_e|)^{j-2}  c_2 (\ve + \delta^\prime)  
                      \cdot 3 |\mathrm{Id}_e|                                
                        + (3 |\mathrm{Id}_e|)^{j-1} \cdot c_2 (\ve + \delta^\prime) 
                          \\
                          &=j (3 |\mathrm{Id}_e|)^{j-1}  c_2 (\ve + \delta^\prime),                                               \end{align*}
which proves \eqref{eq.0121-1}.
Note that we used the assumption of the induction for the third inequality.
                                      
Combining \eqref{eq.0121-2}--\eqref{eq.0121-1}, 
  we finish the proof of the lemma.
 \end{proof}

     
\begin{lemma} \label{lem.gen_dil}
Let $n,m \in \N$, $\tau\in (0,1]$ and 
$A, \hat{A} \in \R^{m \times n}$.
   Let $\omega$ be a control function.
      \\
       {\rm (i)}~If $\mathbf{w} \in G\Omega_p (\R^{n})$ satisfies that
               \[
         |\mathbf{w}_{s,t}^k | \le \omega(s,t)^{k/p}, 
          \qquad 
            0\le s \le t \le  \tau, \,\, 1 \le k \le [p],
        \] 
then we have
  \[
         |(\Gamma_A \mathbf{w})_{s,t}^k | \le   \{ |A| \, \omega(s,t)\}^{k/p}, 
          \qquad 
            0\le s \le t \le  \tau, \,\, 1 \le k \le [p].
        \] 
        \\
        {\rm (ii)}~If $\mathbf{w}, \hat{\mathbf{w}}\in G\Omega_p (\R^{n})$ and $\ve \ge 0$ satisfiy that       
        \begin{eqnarray*}
         |\mathbf{w}_{s,t}^k | \vee  |\hat{\mathbf{w}}_{s,t}^k | \le \omega(s,t)^{k/p}, 
          \,\,
           |\mathbf{w}_{s,t}^k  -\hat{\mathbf{w}}_{s,t}^k |   \le \ve\omega(s,t)^{k/p}, 
            \quad
               0\le s \le t \le  \tau, \,\, 1 \le k \le [p],
        \end{eqnarray*}
        then we have
  \begin{align*}
         |(\Gamma_A \mathbf{w})_{s,t}^k 
            - (\Gamma_{\hat{A}} \hat{\mathbf{w}} )_{s,t}^k| 
            &\le    \{ (|A|\vee |\hat{A}|)^k \ve + k  (|A|\vee |\hat{A}|)^{k-1} |A-\hat{A}|\}
             \, \omega(s,t)^{k/p}
                \end{align*}
           for all $0\le s \le t \le  \tau$ and $1 \le k \le [p]$.
                        \end{lemma}

 \begin{proof}   
 We can easily show this by straightforward computation.
  \end{proof}


\begin{proposition}\label{pr.diff_JK}
Let the notation be  as above.
Then, for every $\beta \in (0,\beta_0]$, there exists $C>0$ such that
\begin{align*}
\lefteqn{
 d_{p\textrm{-}\mathrm{var}} \bigl( (\mathbf{x}, \mathbf{y}, \mathbf{J}, \mathbf{K}),
   (\hat{\mathbf{x}}, \hat{\mathbf{y}}, 
  \hat{\mathbf{J}}, \hat{\mathbf{K}})
     \bigr)
                 }
                  \nn\\ 
                   &\le 
          C(\vertiii{\mathbf{x}}_{p\textrm{-}\mathrm{var}}^p + 
 \vertiii{\hat{\mathbf{x}}}_{p\textrm{-}\mathrm{var}}^p +1)
          \exp \bigl[
            C (N^p_{\beta} (\mathbf{x})+ N^p_{\beta} (\hat{\mathbf{x}}))
            \bigr]
          (d_{p\textrm{-}\mathrm{var}} (\mathbf{x}, \hat{\mathbf{x}}) +|a-\hat{a}|).        
            \end{align*}
            for all $(\mathbf{x}, a),  (\hat{\mathbf{x}}, \hat{a}) 
\in G\Omega_p (\R^d)\times  \R^e$.          Here, the constant $C$ depends only on $\beta$ and $K'$.       
            \end{proposition}    
     
  \begin{proof}   
  Denote by $Q_{A,B} \in L(\R^d \oplus \mathcal{V}, \R^d \oplus \mathcal{V})$
  the linear map introduced in the
   definition of $\tilde{\Gamma}_{A,B}$ in \eqref{def.Gamma_tilde}. 
  It is easy to see that
  \begin{equation}\label{ineq.0123-1}
  |Q_{A,B}| \le  1\vee |A| \vee |B|, \qquad 
        |Q_{A,B}- Q_{\hat{A}, \hat{B}}| \le  |A- \hat{A}| \vee |B- \hat{B}|.
            \end{equation}

       In what follows, we use the same notation 
        as in  the proof of Lemma \ref{lem.diff_1st_Jac}.  
          We write  $R:= \|J\|_{\infty} \vee\|\hat{J}\|_{\infty} \vee \|K\|_{\infty}\vee \|\hat{K}\|_{\infty}$,
          $\delta := \|J -\hat{J}\|_{\infty} \vee \|K-\hat{K}\|_{\infty}$,
          $\delta' := \|y -\hat{y}\|_{\infty}$ 
          and $\ve := d_{p\textrm{-}\mathrm{var}} (\mathbf{x}, \hat{\mathbf{x}})$ for simplicity. (Note that $R \ge 1$.)
          In this proof, $c_i~(i=1,2,\dots)$ are positive constants which depends only on $\beta$ and $K'$.
           
We see from Lemma \eqref{lem.diff_loc_Jac} that 
there exist constants $c_1, c_2 >0$ such that
       \begin{align*}
         |(\mathbf{x}(j), \mathbf{y}(j), \mathbf{J}(j), \mathbf{K}(j))_{s,t}^k |
                  \vee  |(\hat{\mathbf{x}}(j), \hat{\mathbf{y}}(j), 
  \hat{\mathbf{J}}(j), \hat{\mathbf{K}}(j))_{s,t}^k|                        
   &\le c_1\omega(s,t)^{k/p},
                  \\
               |(\mathbf{x}(j), \mathbf{y}(j), \mathbf{J}(j), \mathbf{K}(j))_{s,t}^k 
               -  (\hat{\mathbf{x}}(j), \hat{\mathbf{y}}(j), 
  \hat{\mathbf{J}}(j), \hat{\mathbf{K}}(j))_{s,t}^k| 
               &\le c_2 (\ve +\delta') \omega(s,t)^{k/p}
                     \end{align*}
for all $u_{j-1}\le s \le t \le  u_{j}$, $1 \le j \le M$ and $1 \le k \le [p]$.
      Here, $\omega$ is given by  \eqref{def.diff_omega}.       

Since 
$(\mathbf{x}, \mathbf{y}, \mathbf{J}, \mathbf{K})
=\tilde{\Gamma}_{ J_{\tau_{j-1}}, K_{\tau_{j-1}}}
(\mathbf{x} (j), \mathbf{y}(j), \mathbf{J}(j), \mathbf{K}(j))$,
we can use Lemma \ref{lem.gen_dil} and \eqref{ineq.0123-1} on each subinterval.          
   For certain constants $c_3, c_4>0$,  we have
            \begin{align*}
         |(\mathbf{x}, \mathbf{y}, \mathbf{J}, \mathbf{K})_{s,t}^k |
                  \vee  |(\hat{\mathbf{x}}, \hat{\mathbf{y}}, 
  \hat{\mathbf{J}}, \hat{\mathbf{K}})_{s,t}^k|                        
   &\le \{ c_3R\, \omega(s,t)\}^{k/p},
                  \\
               |(\mathbf{x}, \mathbf{y}, \mathbf{J}, \mathbf{K})_{s,t}^k 
               -  (\hat{\mathbf{x}}, \hat{\mathbf{y}}, 
  \hat{\mathbf{J}}, \hat{\mathbf{K}})_{s,t}^k| 
               &\le  (\ve +\delta' +\delta)\{c_4R\, \omega(s,t)\}^{k/p}
                     \end{align*}
for all $u_{j-1}\le s \le t \le  u_{j}$, $1 \le j \le M$ and $1 \le k \le [p]$.

Finally, we use Lemma \ref{lem.loc_to_glo} to obtain a global estimate:
   For a certain constants $c_5>0$,  it holds that
     \[
     |(\mathbf{x}, \mathbf{y}, \mathbf{J}, \mathbf{K})_{s,t}^k 
               -  (\hat{\mathbf{x}}, \hat{\mathbf{y}}, 
  \hat{\mathbf{J}}, \hat{\mathbf{K}})_{s,t}^k| 
               \le  (\ve +\delta' +\delta) M^{k(p-1)/p} \{c_5 R\, \omega(s,t)\}^{k/p}
                    \]
                    for all $0\le s \le t \le 1$ and $1 \le k \le [p]$. 
                     Note that  $M \le N^p_{\beta} (\mathbf{x})+N^p_{\beta} (\hat{\mathbf{x}}) +2$ and $\omega (0,1) \le \vertiii{\mathbf{x}}_{p\textrm{-}\mathrm{var}}^p + 
 \vertiii{\hat{\mathbf{x}}}_{p\textrm{-}\mathrm{var}}^p  + [p]$.  
  Recall the estimate of $\delta'$, $R$ and $\delta$ in 
  (the proof of)  Corollary \ref{cor.0117-1}, Lemma \ref{lem.glob_est}
    and Lemma \ref{lem.diff_1st_Jac}, respectively.  
         \end{proof}


\section{Fractional Brownian rough path}

We first introduce the abstract Wiener space for 
fractional Brownian motion (fBM) with Hurst parameter $H$.
Throughout this paper we assume $H \in (1/4, 1/2]$.

Let $\mu =\mu^{H}$ be the law of $d$-dimensional fBM with 
Hurst parameter $H$ defined on $\cW := \cC_0 (\R^d)$.
The Cameron-Martin space is denoted by $\cH =\cH^{H}$.
Then, $(\cW, \cH^H, \mu^H)$ becomes an abstract Wiener space. 
(When we need to specify the dimension $d$, we will write 
$(\cW^d, \cH^{H, d}, \mu^{H, d})$.)
The generic element of $\cW$ is denoted by $w$.
By definition, the coordinate process $w=(w_t)_{t\in [0,1]}$
is an fBM with Hurst parameter $H$, that is, 
it is a $d$-dimensional 
mean-zero Gaussian process with the covariance 
\[
\E [w^i_s w^j_t] =\frac12 ( t^{2H} +s^{2H} -|t-s|^{2H}) \delta_{ij},
\qquad\quad 
s,t \in [0,1], \,\, 1 \le i,j \le d,
\]
where $\delta_{ij}$ is Kronecker's delta.
When $H =1/2$, fBM coincides with the usual Brownian motion.
According to \cite[Theorem 1.1]{fggr}, there exists a 
continuous embedding 
\[
\cH^H \hookrightarrow \cC_0^{\rho_{H}\textrm{-}\mathrm{var}} (\R^d)
\qquad
\mbox{with }\,\,
\rho_{H} := \frac{1}{H +(1/2)} \in [1,4/3).
\]
If $p \in (H^{-1}, \infty)$ is close enough to $H^{-1}$, 
then $p^{-1} + \rho_{H}^{-1}>1$ (the condition for Young integration) 
holds.


Since $H\in (1/4,1/2]$, fBM admits a canonical RP lift, 
which is called fractional Brownian RP (fBRP).
fBRP can be obtained as the limit of various approximation methods, 
but in this paper it is viewed as the limit of piecewise linear approximations.

For $w \in \cW$
and the partition $\cP =\{ 0 =t_0 <t_1 <\cdots <t_N= 1\}$,
we denote by $w(\cP) \in\cW$ the piecewise linear approximation 
of $w$ associated with $\cP$, that is, $w_{t_j}= w(\cP)_{t_j}$
for all $0 \le j \le N$ and $w(\cP)$ is linearly
interpolated  on $[t_{j-1}, t_j]$ for all $1 \le j \le N$.
Then, if $H^{-1} <p < [H^{-1}]+1$, 
the family of $G\Omega_p (\R^d)$-valued random 
  variables $\{ S_p (w(\cP)) \}_{\cP}$ converges as $|\cP| \searrow 0$.
    (The limit is denoted by $\mathbf{w}$ and called fBRP
     with Hurst parameter $H$. 
       Note that $\mathbf{w}:= \lim_{|\cP| \searrow 0}S_p (w(\cP))$ is independent of $p$.)  
         More precisely, we have 
         \[
         \lim_{|\cP| \searrow 0}
         \E[ d_{p\textrm{-}\mathrm{var}} (\mathbf{w}, S_p (w(\cP)))^q] =0
                  \]
for every $q \in [1, \infty)$. (See \cite[Subsection 15.5.1]{fvbk} for instance.)

If $H \in (1/3, 1/2]$ (resp. $H \in (1/4, 1/3]$), 
$p$ can be chosen so that $[p]=2$ (resp. $[p]=3$).
This is the most natural choice 
(at least from the viewpoints of topology and work efficiency).
However, as long as we use this kind of roughness $p$, 
  we cannot obtain a satisfactory convergence rate for the piecewise linear 
   approximations.  
     The next proposition, due to Friz-Riedel \cite{fr}, claims that 
       we can obtain a desired convergence rate 
          if we take $p$ sufficiently
            large (i.e. weaken the RP topology sufficiently).

Let $H$ and $p$ as above take any $r \in (p, \infty)$.
Then, $\mathbf{Ext}_{p,r} (\mathbf{w})$ is a $G\Omega_{r} (\R^d)$-valued random variable and is also called fBRP (of roughness $r$).
By the continuity of $\mathbf{Ext}_{p,r}$, 
$\{ S_r (w(\cP)) \}_{\cP}$ converges to $\mathbf{Ext}_{p,r} (\mathbf{w})$
as $|\cP| \searrow 0$ at least in probability.
(In fact, the convergence takes place in every $L^q$, $1 \le q <\infty$.
See the next proposition).
Since no confusion may occur, we will slightly abuse the symbol to
 write 
$\mathbf{w}=(1, \mathbf{w}^1, \ldots, \mathbf{w}^{[r]})$ 
for $\mathbf{Ext}_{p,r} (\mathbf{w})$ for the rest of this paper.

\begin{proposition} \label{pr.FR}
Let the notation be as above.
Then, for every  $\ve \in (0, 2H -\tfrac12)$, there exists $p_\ve >H^{-1}$ 
such that the following inequality holds true:

For every $p \in (p_\ve, \infty)$, there exists a constant $C=C(p, \ve, H, d)>0$
independent of $q$ and $\cP$ such that
\[
\E[ d_{p\textrm{-}\mathrm{var}} (\mathbf{w}, S_p (w(\cP)))^q]^{\frac{1}{q}}
\le 
 C q^{\frac{[p]}{2}} |\cP|^{2H -\frac12 -\ve}
\]
holds for every $q \in [1, \infty)$ and $\cP$.
Here, $\mathbf{w}=(1, \mathbf{w}^1, \ldots, \mathbf{w}^{[p]})$
stands for fBRP of roughness $p$.
\end{proposition}

\begin{proof} 
This is a special case of the main results of \cite{fr}.
Combine Theorem 5 and results in Subsection 6.2 of \cite{fr}.
\end{proof}

\begin{remark} \label{rem.L^q_bound}
The above proposition is borrowed from \cite{fr}
and is a refinement of \cite[Theorem 15.42 {\rm (ii)}]{fvbk},
in which the $L^q$-convergence of  the piecewise linear approximations 
is shown for the roughness $p$ slightly larger than $1/H$.
In the theorem in \cite{fvbk}, the RP topology is idealistic, 
 but  the convergence rate is not the desired one.
The meaning of \cite{fr} is that the desired convergence rate is obtained
at the price of the roughness.

By the way, \cite[Theorem 15.42 {\rm (ii)}]{fvbk} 
(and Lyons' extension) also implies the following 
$L^q$-boundedness of $\{ S_p (w(\cP)) \}_{\cP}$.
For every $p >H^{-1}$ and $q \in [1,\infty)$, 
\[
\E [  \vertiii{\mathbf{w}}_{p\textrm{-}\mathrm{var}}^q]^{\frac{1}{q}}
   +
     \sup_{\cP} 
      \E [  \vertiii{ S_p (w(\cP))}_{p\textrm{-}\mathrm{var}}^q]^{\frac{1}{q}}     <\infty.
                  \]
To see this fact, just consider the case $\cP =\{0,1\}$ (the trivial partition)
 in that theorem.
\end{remark}


For the rest of this section
we will show the exponential integrability of
$N^p_{\beta} (\mathbf{w})$ and $N^p_{\beta} (S_p (w(\cP)))$
for fRBP $\mathbf{w}$, 
 following a well-known result in \cite{cll}.
First, we introduce a condition on  a family of 
$G\Omega_{p} (\R^d)$-valued 
random variables $\{ \mathbf{z}_\alpha\}_{\alpha \in A}$
defined on a probability space $(\Omega, \cF, \mathbb{P})$,
where $A$ is any index set.
\begin{definition}
We say that $\{ \mathbf{z}_\alpha\}_{\alpha \in A}$ satisfies {\bf (ExpN)}
on $G\Omega_{p} (\R^d)$ if we have 
\[ 
\sup_{\alpha \in A}
\E \bigl[ \exp (\eta N^p_{\beta} (\mathbf{z}_\alpha ))
\bigr] <\infty
\qquad
\mbox{for all $\eta >0$ and $\beta >0$.} 
\]
\end{definition}

\begin{example}\label{exm.0117}
Let $H \in (1/4, 1/2]$ and assume that $p\in (H^{-1}, [H^{-1}]+1)$ is 
sufficiently close to $H^{-1}$ (so that $p^{-1} + \rho_{H}^{-1}$ holds).
Then, $\{\mathbf{w}\} \cup \{S_p (w(\cP)) : \cP\}$ is known to satisfy
 {\bf (ExpN)} on $G\Omega_{p} (\R^d)$.
 The exponential integrability of $N^p_{\beta} (\mathbf{w})$ was proved 
 in \cite{cll}.
The  {\bf (ExpN)}-property of $\{S_p (w(\cP)) : \cP\}$ was essentially 
proved in \cite[Section 5]{ina14}. 
(The difference is that only the dyadic partitions are considered in \cite{ina14}.
But, this is irrelavent and the same proof works for general $\cP$.)
 \end{example}


\begin{lemma}\label{lem.ExpN1}
Let $2 \le p <4$ and
assume that $\{ \mathbf{z}_\alpha\}_{\alpha \in A}$ satisfies {\bf (ExpN)}
on $G\Omega_{p} (\R^d)$.
Then, $\{ (\mathbf{z}_\alpha,  \bm{\lambda})\}_{\alpha \in A}$ satisfies {\bf (ExpN)}
on $G\Omega_{p} (\R^{d+1})$.
Here, $(\mathbf{z}_\alpha,  \bm{\lambda})$ is the Young pairing 
of $\mathbf{z}_\alpha$ and $\lambda$. (Recall that $\lambda_t =t$.) 
\end{lemma}

\begin{proof}
By straightforward computation, we have 
\[
\vertiii{\mathbf{(\mathbf{x},  \bm{\lambda})}}_{p\textrm{-}\mathrm{var}, [s,t]}
\le 
c(\vertiii{\mathbf{x}}_{p\textrm{-}\mathrm{var}, [s,t]}
+
\|\lambda\|_{1\textrm{-}\mathrm{var}, [s,t]})
=
c\{
\vertiii{\mathbf{x}}_{p\textrm{-}\mathrm{var}, [s,t]}
+ (t-s )\}
\]
for some constant $c=c(p)>0$ independent of 
$\mathbf{x}\in G\Omega_{p} (\R^d)$ and $(s,t)\in\triangle$,
which may vary from line to line.
   It immediately follows that 
    \[
     \vertiii{\mathbf{(\mathbf{x},  \bm{\lambda})}}_{p\textrm{-}\mathrm{var}, [s,t]}^p
\le c \{   \vertiii{\mathbf{x}}_{p\textrm{-}\mathrm{var}, [s,t]}^p\vee (t-s)\}. 
  \] 
  This implies that 
 $N^p_{\beta} ( (\mathbf{x},  \bm{\lambda})) 
  \le N^p_{\beta/c} (\mathbf{x}) +[c/\beta]+1$.  
    This proves the lemma.
\end{proof}

\begin{lemma}\label{lem.ExpN2}
Assume that $\{ \mathbf{z}_\alpha\}_{\alpha \in A}$ satisfies {\bf (ExpN)}
on $G\Omega_{p} (\R^d)$ for $p\in [2, \infty)$.
Then, $\{\mathbf{z}_\alpha\}_{\alpha \in A}$ satisfies {\bf (ExpN)}
on $G\Omega_{r} (\R^d)$ for every $r \in (p,\infty)$.
\end{lemma}

\begin{proof}
Let $\mathbf{x}\in G\Omega_{p} (\R^d)$ be arbitrary and
write $\omega (s,t) = \vertiii{\mathbf{x}}_{p\textrm{-}\mathrm{var}, [s,t]}^p$.
It is easy to see that, for all $1 \le k\le [p]$ and $(s,t)\in\triangle$, 
  we have $\|\mathbf{x}\|_{r/k\textrm{-}\mathrm{var}, [s,t]}^{r/k}
\le \omega (s,t)^{r/p}$.
By Lyons' extension theorem (see Inequality \eqref{eq.0113-1}),
we also have 
\[
 \|\mathbf{x}\|_{r/k\textrm{-}\mathrm{var}, [s,t]}^{r/k}
 \le
\{\beta_p  (k/p)!\}^{-r/k} \omega (s,t)^{r/p}.
\]
for $[p]+1\le k \le [r]$ and $(s,t)\in\triangle$.
So, there exists a constant $c =c(p,r)>0$ such that 
\[
\vertiii{\mathbf{x}}_{r\textrm{-}\mathrm{var}, [s,t]}^r
  \le c (\vertiii{\mathbf{x}}_{p\textrm{-}\mathrm{var}, [s,t]}^p)^{r/p}, 
  \qquad (s,t)\in\triangle.
    \]
Therefore,  if $\vertiii{\mathbf{x}}_{r\textrm{-}\mathrm{var}, [s,t]}^r >\beta$, 
then $\vertiii{\mathbf{x}}_{p\textrm{-}\mathrm{var}, [s,t]}^p >\beta' := (\beta /c)^{r/p}$.
This implies that $N^r_{\beta} (\mathbf{x}) \le N^p_{\beta'} (\mathbf{x})$.
Our assertion immediately follows from this.
\end{proof}

\begin{proposition}\label{pr.fbrp_ExpN}
For $H \in (1/4, 1/2]$ and $p\in (H^{-1}, \infty)$, we have the following:
\\
{\rm (i)} $\{\mathbf{w}\} \cup \{S_p (w(\cP)) : \cP\}$  satisfies
 {\bf (ExpN)} on $G\Omega_{p} (\R^d)$ 
 \\
 {\rm (ii)} $\{(\mathbf{w}, \bm{\lambda})\} \cup \{(S_p (w(\cP)), \bm{\lambda}) : \cP\}$  satisfies
 {\bf (ExpN)} on $G\Omega_{p} (\R^{d+1})$. 
 \end{proposition} 

\begin{proof}
Just combine Example \ref{exm.0117}
and Lemmas \ref{lem.ExpN1} and \ref{lem.ExpN2}.
(It should also be noted that taking the Young pairing with $\lambda$ 
and Lyons' extension map commute.)
\end{proof}

\section{ODE driven by piecewise linear approximation}
\label{sec.pla}

In this section, we consider ODEs driven by a piecewise linear 
approximation of $w$.
In particular, derivatives of the solution maps of the ODEs 
are studied.
We basically follow  arguments in \cite[Section 3]{ina14}.
All line integrals and ODEs in this section 
are  understood  in the Riemann-Stieltjes sense.
From now on we will assume that 
$\sigma \in C_{\mathrm{b}}^{\infty} (\R^e, \R^{e\times d})$
and $b\in C_{\mathrm{b}}^{\infty} (\R^e, \R^e)$.
We will also
 write $\sigma :=[ V_1,\ldots, V_d] \colon \R^e \to \R^{e\times d}$ and $b =V_0$.

We consider the following ODE on $\R^e$
driven by $w\in\cC_0^{1\textrm{-}\mathrm{var}} (\R^d)$ which corresponds to RDE \eqref{def.introRDE} (with $a=0$):
\begin{equation}
dy_t = \sigma (y_t) dw_t + b  (y_t)  dt, 
\qquad
\quad  y_0 =0.
\label{def.RSode8}
\end{equation} 
Or equivalently, 
\begin{equation} \label{def.RSode1}
dy_t = \sum_{j=1}^d V_j (y_t) dw_t^j + V_0 (y_t) dt, 
\qquad
\quad  y_0 =0.
\end{equation}
We write $y_t =I_t (w)$. 
  It is well-known  that $I\colon \cC_0^{1\textrm{-}\mathrm{var}} (\R^d)\to \cC_0^{1\textrm{-}\mathrm{var}} (\R^e)$ is Fr\'echet smooth.

The associated Jacobian ODE and its inverse are given as follows:
\begin{align} 
dJ_t &= \sum_{j=1}^d \nabla V_j (y_t) J_t dw_t^j  + \nabla V_0 (y_t) J_t dt,
\qquad
\quad  J_0 =\mathrm{Id}_e,
\label{def.RSode2} 
\\
dK_t &= -  \sum_{j=1}^d K_t \nabla V_j (y_t) dw_t^j -K_t \nabla V_0 (y_t)dt,
\qquad
\quad  K_0 =\mathrm{Id}_e.
\label{def.RSode3} 
\end{align} 
These are $\R^{e\times e}$-valued ODEs.
Here, $\nabla V_j$ is viewed as an $\R^{e \times e}$-valued function. 
Note that $J_t = K_t^{-1}$ always holds.
Since a Gronwall-type lemma holds in $1$-variational setting, 
the system of ODEs \eqref{def.RSode1}--\eqref{def.RSode3}
has a unique (global) solution and, moreover, the inequality
$\|J\|_{1\textrm{-}\mathrm{var}} + \|K\|_{1\textrm{-}\mathrm{var}} \le C \exp(C\|w\|_{1\textrm{-}\mathrm{var}})$ holds for a constant $C>0$ independent of $w$. 


\begin{remark}   
Here, we prefer writing the equations for $(y, J, K)$ in essentially
the same way as in Remark \ref{rem.equivRDE}.
Of course, writing them in an analogous way to 
\eqref{def.JRDE_1}--\eqref{def.JRDE_3} is possible and equivalent. 
\end{remark}

\begin{remark}\label{rem.ode_double}
The system of ODEs \eqref{def.RSode1} --\eqref{def.RSode3} 
 is driven by $w$. 
  However, one can also view that the system is driven by 
   $(w, \theta) \in \cC_0^{1\textrm{-}\mathrm{var}} (\R^{d+d})=\cC_0^{1\textrm{-}\mathrm{var}} (\R^d)^{\oplus 2}$. 
    (The coefficient vector fields in front of $d\theta^j$'s are identically zero.)
    \end{remark}


Now we write down the directional derivatives of $I$
in the direction of  $h \in \cC_0^{1\textrm{-}\mathrm{var}} (\R^d)$. 
It is known that
    $D_h^n y_t =(D_h)^n y_t = D^n I_t (w)\la h, \ldots, h\ra$ satisfies 
      the following ODEs for all $n\in\N$.
(Needless to say, this differentiation $D$ is with respect to the $w$-variable.)

\begin{align} 
D_h y_t &=
   J_t \int_0^t K_s \sigma (y_s) dh_s, 
\label{eq.0131-1}
\\
D_h^2 y_t &=
   J_t \int_0^t K_s \{ \nabla^2\sigma (y_s) \la D_h y_s, D_h y_s, dw_s \ra
    +  2 \nabla \sigma (y_s) \la D_h y_s, dh_s \ra
     \nn\\
     &\qquad\qquad\qquad
               +  \nabla^2 b(y_s) \la D_h y_s, D_h y_s \ra ds\}.
   \label{eq.0131-2}
    \end{align}
Here, $\nabla$ stands for the standard gradient operator on $\R^e$.
For general $n \ge 2$, 
\begin{align} 
D_h^n y_t &=
 J_t \int_0^t K_s 
\Bigl\{
\sum_{l=2}^n   
\sum_{
i_1 + \ldots + i_l =n}
C_{i_1, \ldots, i_l}
\nabla^l \sigma ( y_s) \la D_h^{i_1} y_s ,  \ldots , D_h^{i_l} y_s, dw_s\ra  
\nn\\
&
\qquad\qquad\quad  +
\sum_{l=1}^{n-1}   
\sum_{
i_1 + \ldots + i_l =n-1}
C'_{i_1, \ldots, i_l}
\nabla^l \sigma ( y_s) \la D_h^{i_1} y_s ,  \ldots , D_h^{i_l} y_s, dh_s\ra  
\nn\\
&
\qquad\qquad \quad +
\sum_{l=2}^n   
\sum_{
i_1 + \ldots + i_l =n}
C_{i_1, \ldots, i_l}
\nabla^l b ( y_s) \la D_h^{i_1} y_s ,  \ldots , D_h^{i_l} y_s\ra ds
\Bigr\}.
   \label{eq.0131-3}
    \end{align}
      Here, {\rm (i)}~ the summation $\sum_{i_1 + \ldots + i_l =n}$
runs over all non-decreasing sequence $0< i_1 \le  \ldots \le i_l$ of natural numbers
such that $i_1 + \ldots + i_l =n$,
{\rm (ii)}~
$C_{i_1, \ldots, i_l},  C'_{i_1, \ldots, i_l} \in {\mathbb N}$
are constants, 
but their exact values are not important for our purpose.

\begin{remark}\label{rem.0203-1}
In what follows, we will sometimes write $\Xi_{n,t} (w,h)
$ for $D_h^n y_t (w)$, $n\in\N$.
As one can easily see, $(w, h) \mapsto \Xi_{n, \cdot} (w,h)$
is locally Lipschitz continuous from
$\cC_0^{1\textrm{-}\mathrm{var}} (\R^{2d})$ to $\cC_0^{1\textrm{-}\mathrm{var}} (\R^{e})$
 for every $n$.
     For a partition $\cP$, we write 
       $\Xi_n (\cP)_t (w, h) = \Xi_{n, t} (w(\cP),h(\cP))$. 
     \end{remark}

\vspace{10mm}
    
 For the canonical realization $w=(w_t)_{t \in [0,1]}$ of $d$-dimensional 
 fBM with $H\in (1/4, 1/2]$
and the partition $\cP =\{ 0 =t_0 <t_1 <\cdots <t_N= 1\}$,    
we replace the driver by the piecewise linear approximation 
$w(\cP)=(w(\cP)_t)_{t \in [0,1]}$ of fBM.   
    The solution of  ODEs \eqref{def.RSode1}--\eqref{def.RSode3}
     driven by $w(\cP)$ is denoted by $(y(\cP), J(\cP), K(\cP))$.    
   Those are Wiener functionals on the 
   abstract Wiener space $(\cW, \cH^H, \mu^H)$.    
    
\begin{lemma} \label{lem.aprox_dinfty}
For every $\cP$ and $t \in [0,1]$, we have $y(\cP)_t \in \mathbf{D}_{\infty} (\R^e)$.
\end{lemma}
    
\begin{proof} 
This lemma was essentially proved in \cite[Proposition 3.2]{ina14}.
The only differences are that (1) the partition $\cP$ is arbitrary
(not just the dyadic ones)
and (2) the ODE for $y$ has a drift term.
But, the same proof works anyway. So, we only give a sketch of proof here.

From the Fr\'echet smoothness of $I$, we see that 
\[
y(\cP)_t  = G \bigl(\{ w_{t_j} - w_{t_{j-1}} \}_{j=1}^N  \bigr)
\]
holds for a certain $G \in C^\infty (\R^{dN}, \R^e)$.
Since 
\[
\|J(\cP)\|_{1\textrm{-}\mathrm{var}} + \|K(\cP)\|_{1\textrm{-}\mathrm{var}}
 \le C \exp \bigl(C\sum_{j=1}^N | w_{t_j} - w_{t_{j-1}}| \bigr)
 \]
  holds for some constant $C=C(\cP) >0$, we can show that 
   $\nabla^k G$ is at most of exponential growth
   ($k \in \N \cup \{0\}$).  
    From \cite[Corollary 5.3.2]{mt}, we can see that 
   $y(\cP)_t \in \mathbf{D}_{\infty} (\R^e)$.
(This last part is heuristically obvious, but is not completely obvious.
A rigorous proof can be found in the above mentioned reference.)
\end{proof}


Let 
$(\theta_t) = (\theta_t^1, \ldots, \theta_t^d)_{t \in [0,1]}$ be an independent copy of 
$(w_t)$.
The abstract Wiener space that corresponds to 
the $2d$-dimensional fBM $(w_t, \theta_t)_{t\in [0,1]}$ is 
\[
(\cW^{2d}, \cH^{H, 2d}, \mu^{H, 2d}) =
(({\cal W}^d)^{\oplus 2}, ({\cal H}^{H,d})^{\oplus 2}, \mu^{H,d} \times \mu^{H,d}).
\]
The expectation with respect to $w$-variable and $\theta$-variable
are
denoted by ${\mathbb E}'$ and $\hat{\mathbb E}$, respectively.
The expectation with respect to $(w, \theta)$-variable 
is denoted by
${\mathbb E}= {\mathbb E}^{\mu \times \mu} = {\mathbb E}' \times \hat{\mathbb E}$.

If $w$ is shifted by $h$, then $w(\cP)$ is shifted by $h(\cP)$. 
So, we have 
\begin{equation}
\label{eq.0201-1}
D^n y(\cP)_t \la  h,h, \ldots, h\ra 
=
D^n I_t (w(\cP)) \la  h(\cP), h(\cP), \ldots, h(\cP)\ra
=\Xi_n (\cP)_t (w, h),
\end{equation}
where $D$ in front of $I_t (w(\cP))$ is in the Fr\'echet sense, 
while 
$D$   in front of $y(\cP)_t$
is ${\cal H}$-derivative in the sense of Malliavin calculus.


In \eqref{eq.0201-1}, $h$ is still deterministic.
By replacing $h(\cP)$ in \eqref{eq.0201-1} with 
  random piecewise linear path $\theta(\cP)$, we have
\begin{equation}\label{eq.0201-2}
\Xi_n (\cP)_t (w,\theta) 
=
D^n I_t (w(\cP)) \la  \theta(\cP), \theta(\cP), \ldots, \theta(\cP)\ra,
\end{equation}
which is a cylinder functional of the $2d$-fBM $(w_t, \theta_t)_{t \in [0,1]}$.
More explicitly, 
\begin{equation}
 \Xi_1 (\cP)_t
=
J(\cP)_t \int_0^t K(\cP)_s  \sigma ( y(\cP)_s) d\theta(\cP)_s
\label{eq.0201-3}
\end{equation}
and, for $n\ge 2$,
\begin{align}
\Xi_n (\cP)_t
&= 
J(\cP)_t \int_0^t K(\cP)_s
\nn\\
&\quad 
\times
\Bigl\{
\sum_{l=2}^n   
\sum_{
i_1 + \ldots + i_l =n}
C_{i_1, \ldots, i_l}
\nabla^l \sigma ( y(\cP)_s) \la  \Xi_{i_1} (\cP)_s ,  \ldots ,  \Xi_{i_l} (\cP)_s, dw(\cP)_s\ra  
\nn\\
&
\qquad+
\sum_{l=1}^{n-1}   
\sum_{
i_1 + \ldots + i_l =n-1}
C'_{i_1, \ldots, i_l}
\nabla^l \sigma ( y(\cP)_s) \la  \Xi_{i_1} (\cP)_s,  \ldots , \Xi_{i_l} (\cP)_s, d\theta (\cP)_s\ra  
\nn\\
&
\qquad+
\sum_{l=2}^n   
\sum_{
i_1 + \ldots + i_l =n}
C_{i_1, \ldots, i_l}
\nabla^l b ( y(\cP)_s) \la  \Xi_{i_1} (\cP)_s ,  \ldots ,  \Xi_{i_l} (\cP)_s\ra ds  
\Bigr\},
\label{eq.0201-4}
\end{align}
where the positive constants are the same as in (\ref{eq.0131-3}).

Since $\Xi_n (\cP)_t (w,\theta)$  is defined for all $w$ and $\theta$,
we can think of $\Xi_n (\cP)_t (w,\,\cdot\,)$ as a Wiener functional in $\theta$ for each fixed $w$.
Then, it is clear from the right hand side of (\ref{eq.0201-2})
that, for each $w$,  $\Xi_n (\cP)_t (w,\theta)$ is a polynomial of order $n$
in $\{\theta_{t_i} : t_i \in \cP\}$.
In particular,  
$\Xi_n (\cP)_t (w,\,\cdot\,)$ belongs to $n$th order inhomogeneous Wiener chaos.
Moreover, by straight-forward computation, we have
\begin{equation}\label{eq.0202-1}
\hat{D}^n \Xi_n (\cP)_t (w,\,\cdot\,) \la h, h, \ldots, h \ra
=
n! D^n y(\cP)_t \la  h,h, \ldots, h\ra, 
\qquad
h \in {\cal H}^{H, d}.
\end{equation}
Here, $\hat{D}$ stands for ${\cal H}$-derivative with respect to the $\theta$-variable.
Note that neither side depends on $\theta$.
As we remarked in Introduction, 
since both are symmetric $n$-multilinear functional, it holds that 
$\hat{D}^n \Xi_n (\cP)_t (w,\,\cdot\,) = n! D^n y(\cP)_t 
 \in \{({\cal H}^{H, d})^*\}^{\otimes n}$ for each $w$.


    The next proposition shows that, if $\{ \Xi_n(\cP)_t \}_\cP$ is Cauchy
    as $|\cP| \searrow 0$ in $L^q$-norm, 
then $\{D^n y(\cP)_t\}_\cP$ is also Cauchy  in $L^q$-norm.
Moreover,  the convergence rate of $\{ \Xi_n(\cP)_t \}_\cP$
in $L^q$-norm dominates that of $\{D^ny(\cP)_t\}_\cP$ in $L^q$-norm.
Consequently, in order to obtain the convergence rate of $\{y(\cP)_t\}_\cP$ in 
an arbitrary Sobolev norm, 
we only need the convergence rate of
 $\{ \Xi_n(\cP)_t \}$ in $L^q$-norm for every $t, n$ and  $q \in [2, \infty)$.
\begin{proposition}\label{pr.est_ym}
For every $n\in\N$ and $q \in [2, \infty)$,
there is a constant $C=C_{q,n} >0$ independent of $t$ and $\cP$
such that
\[
{\mathbb E}'  [ \|D^n y(\cP)_t  \|_{{\cal H}^{* \otimes n} \otimes {\mathbb R}^e}^q]^{1/q}
\le 
C  {\mathbb E}  [  | \Xi_n(\cP)_t |^q]^{1/q}
\]
for all $0 \le t \le 1$ and $\cP$. 
Here, we wrote $\cH = \cH^{H,d}$ for brevity.
In essentially the same way, we also have 
\[
{\mathbb E}'  [ \|D^n y(\cP)_t - D^n y(\hat\cP)_t  \|_{{\cal H}^{* \otimes n}\otimes {\mathbb R}^e}^q]^{1/q}
\le 
C  {\mathbb E}  [  | \Xi_n(\cP)_t - \Xi_n(\hat\cP)_t |^q]^{1/q}
\]
for all $0 \le t \le 1$ and $\cP, \hat\cP$.
\end{proposition}

\begin{proof}
We prove the first assertion.
From \eqref{eq.0202-1} and Lemma \ref{lem.L_sym},  we see that 
\[
\|D^n y(\cP)_t  \|_{{\cal H}^{* \otimes n} \otimes {\mathbb R}^e} 
= 
\frac{1}{n!} \| \hat{D}^n \Xi_n (\cP)_t (w,\,\cdot\,) \|_{{\cal H}^{* \otimes n} \otimes {\mathbb R}^e} 
 \le
 C_{q,n} \hat{\mathbb E} [  | \Xi_n (\cP)_t (w,\,\cdot\,) |^2]^{\frac12}
   \]
  for some constant $C_{q,n} >0$.
Here, we used the fact that
$\hat{D}^n \Xi_n (\cP)_t (w,\,\cdot\,)$ does not depend on $\theta$.
Taking $L^q$-norm of this inequality with respect to $w$-variable
and using $q/2 \ge 1$, 
we show the first assertion. 

We can show the second assertion 
in the same way, too.
\end{proof}

\section{Wong-Zakai approximation in Sobolev norms}    
    
In this section we show that the Wong-Zakai approximation 
holds (at a fixed time) 
at a desired convergence rate in every Sobolev spaces.
Our main result is immediate from this 
and Hu-Watanabe's approximation theorem (Theorem \ref{thm.hw2.1}).

As in the previous section,
let $(({\cal W}^d)^{\oplus 2}, ({\cal H}^{H,d})^{\oplus 2}, \mu^{H,d} \times \mu^{H,d})$ 
be the abstract Wiener space for the $2d$-dimensional fBM
with Hurst parameter $H \in (1/4, 1/2]$ 
(in the product form).
The coordinate process is denoted by $(w,\theta)= (w_t, \theta_t)_{t\in [0,1]}$.
Then, $(w_t)$ and $(\theta_t)$ are two independent copies 
of the $d$-dimensional fBM with Hurst parameter $H$. 
\begin{definition}\label{def.PLA}
Let $H\in (1/4, 1/2]$, $p > H^{-1}$ and let $F \colon G\Omega_p (\R^{2d})
\to G\Omega_p (\R^m)$ 
be a continuous map for some $m\in\N$. 
We say that $F$ satisfies the {\it piecewise linear approximation 
property} ({\bf (PLA)} for short) with the rate $\delta >0$ if the following 
statement holds:
For every $q \in [1, \infty)$, there exists a constant $C_q >0$
such that
\begin{align} 
\E [ \vertiii{F ( S_p (w(\cP), \theta (\cP)) )  }^q_{p\textrm{-}\mathrm{var}} ]^{1/q} &\le C_q,
\label{def.0203-1}
\\
 \E [ d_{p\textrm{-}\mathrm{var}} \bigl(F ( S_p (w(\cP), \theta (\cP)) ), \,  
   F ( S_p (w(\hat{\cP}), \theta (\hat{\cP})) ) \bigr)^q ]^{1/q} 
 &\le C_q (|\cP|\vee |\hat\cP|)^\delta.
\label{def.0203-2}
\end{align}
hold for all $\cP$ and $\hat\cP$.
Here,  $C_q =C_{q,p, F,\delta}>0$ is a constant independent of $\cP$ and $\hat\cP$. 
\end{definition}
    
\begin{remark} 
Suppose that {\bf (PLA)} holds with the rate $\delta >0$ as in Definition \ref{def.PLA}.
   As we commented in Remark \ref{rem.L^q_bound}, 
         $S_p (w(\hat{\cP}), \theta (\hat{\cP}))$ converges to 
           $(\mathbf{w}, \bm{\theta})$ in any $L^q$-norm, $q \in [1, \infty)$.
            (Here, $(\mathbf{w}, \bm{\theta})$ is viewed as 
              random RPs of roughness $p$ via Lyons' extension.)
Then,  Fatou's lemma and the continuity of $F$ imply that, for all $q \in [1, \infty)$ and $\cP$, 
\begin{align} 
\E [ \vertiii{F (\mathbf{w}, \bm{\theta} ) }^q_{p\textrm{-}\mathrm{var}} ]^{1/q} &\le C_q,
\nn
\\
 \E [ d_{p\textrm{-}\mathrm{var}} \bigl(
 F ( \mathbf{w}, \bm{\theta}), \,\,  
   F ( S_p (w(\cP), \theta (\cP)) )
     \bigr)^q ]^{1/q} 
 &\le C_q |\cP|^\delta.
\nn
\end{align}
Here,  the value of $C_q =C_{q,p, F, \delta}>0$ may be different from that in Definition \ref{def.PLA}. 
\end{remark}

\begin{example} \label{exm.0203-1}
Here are two simple examples of $F$ with {\bf (PLA)} property.
For a given $\ve \in (0, 2H -\tfrac12)$, let $p \ge p_{\ve} >H^{-1}$  
    as in Proposition \ref{pr.FR}.  
     (In this example, by slightly abusing the notation, we denote by 
     $(\mathbf{w}, \bm{\theta} )$  a generic element 
          of $G\Omega_{p} (\R^{2d})$.
      As before, $\lambda_t =t$.)
    \\
(1)~We see from
     Proposition \ref{pr.FR} and Remark \ref{rem.L^q_bound} that
      the identity map of $G\Omega_{p} (\R^{2d})$ 
      (that is, $(\mathbf{w}, \bm{\theta} ) \mapsto (\mathbf{w}, \bm{\theta} )$) satisfies
        {\bf (PLA)} with the rate $2H -(1/2) -\ve$.
\\
(2)~The Young pairing with $\lambda$,
  that is, \[
  (\mathbf{w}, \bm{\theta} ) \in G\Omega_{p} (\R^{2d})\mapsto (\mathbf{w}, \bm{\theta},  \bm{\lambda})\in G\Omega_{p} (\R^{2d+1})
  \]
   satisfies
        {\bf (PLA)} with the rate $2H -(1/2) -\ve$. 
         This fact can be easily verified from (1) since this map is quite  simple. 
                \end{example}


Consider the system of  Riemann-Stieltjes ODEs 
\eqref{def.RSode1}--\eqref{def.RSode3}.
We  view the system is 
driven not just by $w$, but 
by $(w, \theta)\in\cC_0^{1\textrm{-}\mathrm{var}} (\R^{2d})$
(see Remark \ref{rem.ode_double}).
Then the map 
\[
\cC_0^{1\textrm{-}\mathrm{var}} (\R^{2d}) \ni (w, \theta)
\mapsto 
 (w,\theta, \lambda, y, J_\cdot -\mathrm{Id}_e, K_\cdot -\mathrm{Id}_e)
     \in \cC_0^{1\textrm{-}\mathrm{var}} (\R^{2d+1} \oplus\mathcal{V}),
     \]
  where $\lambda_t =t$ and $\mathcal{V}:=\R^e \oplus \R^{e \times e} 
\oplus \R^{e \times e}$,  
    uniquely extends to a continuous map between geometric RP spaces 
     for any $p \in [2, \infty)$:
   \begin{equation}\label{eq.0205-1}
G\Omega_{p} (\R^{2d}) \ni (\mathbf{w}, \bm{\theta} ) 
\mapsto (\mathbf{w}, \bm{\theta}, \bm{\lambda}; \mathbf{y}, \mathbf{J}, \mathbf{K})
\in G\Omega_{p} (\R^{2d+1} \oplus \mathcal{V}).
\end{equation}
  In fact, this continuous map is the Lyons-It\^o map of the corresponding 
  system of RDEs with the initial value 
       $(0, \mathrm{Id}_e, \mathrm{Id}_e)$.
    (See Remark \ref{rem.equivRDE}.
        More precisely, the driving RP $\mathbf{x}$ is replaced by $(\mathbf{w}, \bm{\theta}, \bm{\lambda})$. The coefficient vector field 
    of $dx^j$ is $V_j$ for all $j$,  that of $dt$ is $V_0$, and that of $d\theta^j$ is identically $0$ for all $j$.)
    
\begin{lemma} \label{lem.0203-1}
For a given $\ve \in (0, 2H -\tfrac12)$, let $p \ge p_{\ve} >H^{-1}$  
    as in Proposition \ref{pr.FR}.  
     We use the same notation as in Example \ref{exm.0203-1}.
             Then, the map \eqref{eq.0205-1}
satisfies {\bf (PLA)} with the rate $2H -(1/2) -\ve$. 
\end{lemma}
   
\begin{proof} 
We only need to
 combine Example \ref{exm.0203-1}, Proposition \ref{pr.fbrp_ExpN},
Lemma \ref{lem.glob_est} and Proposition \ref{pr.diff_JK}.
\end{proof}


By the way, $\Xi_n =\Xi_n (w,\theta)$ is defined in 
\eqref{eq.0131-1}--\eqref{eq.0131-3} and Remark \ref{rem.0203-1},
the map 
\begin{align*}
\cC_0^{1\textrm{-}\mathrm{var}} (\R^{2d}) &\ni (w, \theta)
\mapsto 
 \\
 &(w,\theta, \lambda, y, J_\cdot -\mathrm{Id}_e, K_\cdot -\mathrm{Id}_e,
   \Xi_1, \ldots, \Xi_n)
     \in \cC_0^{1\textrm{-}\mathrm{var}} (\R^{2d+1} \oplus\mathcal{V} \oplus (\R^e)^{\oplus n})
     \end{align*}
uniquely extends to a continuous map from $G\Omega_{p} (\R^{2d}) $ 
to $G\Omega_{p} (\R^{2d+1} \oplus \mathcal{V} \oplus (\R^e)^{\oplus n} )$
for any $p \in [2, \infty)$.  We denote the correspondence by 
$(\mathbf{w}, \bm{\theta} ) 
\mapsto (\mathbf{w}, \bm{\theta}, \bm{\lambda}; \mathbf{y}, 
\mathbf{J}, \mathbf{K}, \bm{\Xi}_1, \ldots, \bm{\Xi}_n)$.
  To see this, observe that $\bm{\Xi}_n$ is obtained as an RP integral of 
  $(\mathbf{w}, \bm{\theta}, \bm{\lambda}; \mathbf{y}, 
\mathbf{J}, \mathbf{K}, \bm{\Xi}_1, \ldots, \bm{\Xi}_{n-1})$.
(See also Remarks \ref{rem.RPint1} and \ref{rem.RPint2}).

\begin{proposition} \label{pr.0205-1}
For a given $\ve \in (0, 2H -\tfrac12)$, let $p \ge p_{\ve} >H^{-1}$  
    as in Proposition \ref{pr.FR}.  
     We use the same notation as in Example \ref{exm.0203-1}.
             Then, for every $n\in\N$, the continuous map
\[
G\Omega_{p} (\R^{2d}) \ni (\mathbf{w}, \bm{\theta} ) 
\mapsto (\mathbf{w}, \bm{\theta}, \bm{\lambda}; \mathbf{y}, 
\mathbf{J}, \mathbf{K}, \bm{\Xi}_1, \ldots, \bm{\Xi}_n)
\in G\Omega_{p} (\R^{2d+1} \oplus \mathcal{V} \oplus (\R^e)^{\oplus n} )
\] 
satisfies {\bf (PLA)} with the rate $2H -(1/2) -\ve$. 
Here, we set $\mathcal{V}:=\R^e \oplus \R^{e \times e} 
\oplus \R^{e \times e}$ for simplicity.
\end{proposition}

\begin{proof}
We prove by mathematical induction.  
The case $n=0$ is done in Lemma \ref{lem.0203-1} above.
Suppose that the assertion holds for $n-1$. 
Then, $(\mathbf{w}, \ldots, \bm{\Xi}_n)$ is obtained as an RP 
integral of $(\mathbf{w}, \ldots, \bm{\Xi}_{n-1})$.
Moreover, the integrand and its derivatives are at most of polynomial growth.
(Precisely, the integrand satisfies the assumption of Proposition \ref{pr.RPint} {\rm (ii)}.)
From that proposition and Remark \ref{rem.RPint1}, we see that the 
 the assertion holds for $n$, too.  This completes the proof.
\end{proof}


In Theorems \ref{thm.WZdinfty} and \ref{thm.main},
we denote by $(y_t = \mathbf{y}^{1}_{0,t})_{t \in [0,1]}$
 the first level path of the  solution of RDE \eqref{def.introRDE} (with $a=0$)
driven by fBRP $\mathbf{w}$ of Hurst parameter $H \in (1/4,1/2]$.
Of course, $(y_t)$ does not depend on the choice of 
roughness parameter $p>H^{-1}$ for the RDE.

\begin{theorem}\label{thm.WZdinfty}
For $H \in (1/4,1/2]$, let the notation and situation be as above. 
   Then, the following two statements are satisfied for every 
   $\kappa \in (0, 2H -\tfrac12)$ and $q \in (1, \infty)$:
\\
{\rm (i)}~There exists a constant $C=C_{q,\kappa} >0$ independent of $\cP$ 
such that 
\[
\bigl\| \|y - y(\cP)\|_\infty \bigr\|_{L^q} \le C |\cP|^\kappa
\]
holds for all $\cP$. 
Here, $L^q$-norm is with respect to the law of fBM $(w_t)$. 
\\
{\rm (ii)}~For every $n\in \N$,
There exists a constant $C'=C'_{q,n,\kappa} >0$ independent of 
$\cP$ and $t$ such that 
\[
\| y_t -   y(\cP)_t\|_{\mathbf{D}_{q,n}} \le C' |\cP|^\kappa
\]
holds for all $\cP$ and $0\le t \le 1$. 
Here, the norm is that of $\mathbf{D}_{q,n} (\R^e)$ defined on 
the abstract Wiener space of fBM $(w_t)$.  In particular, $y_t \in \mathbf{D}_\infty (\R^e)$ for all $0\le t \le 1$.
\end{theorem}
    
\begin{proof}  
Define $\ve$ so that $\kappa = 2H -(1/2)-\ve$ and let $p_\ve > H^{-1}$ 
as in Proposition \ref{pr.0205-1}.    
Then, {\rm (i)} is already obtained since the sup-norm is weaker than any 
$p$-variation norm.

We show {\rm (ii)} by estimating the $L^q$-norm of $D^n y(\cP)_t$.
 By Propositions \ref{pr.est_ym} and \ref{pr.0205-1}, we have 
\[
{\mathbb E}'  [ \|D^n y(\cP)_t - D^n y(\hat\cP)_t  \|_{{\cal H}^{* \otimes n}\otimes {\mathbb R}^e}^q]^{1/q}
\le 
C_1  {\mathbb E}  [  | \Xi_n(\cP)_t - \Xi_n(\hat\cP)_t |^q]^{1/q}
\le C_2 (|\cP|\vee |\hat\cP|)^\kappa
\]
for certain constants $C_1, C_2 >0$ independent of $\cP$ and $t$.
Here, we wrote $\cH = \cH^{H,d}$ for brevity.
Hence, $\{D^n y(\cP)_t \}_\cP$ is Cauchy in $L^p$ as $|\cP| \searrow 0$.
 Since $y(\cP)_t \to y_t$ in $L^q$ as $|\cP| \searrow 0$
   and $D$ is closable,  we see that $D^n y_t$ exists 
   and $D^n y(\cP)_t \to D^n y_t$  in $L^q$ as $|\cP| \searrow 0$.  
     By Fatou's lemma and Meyer's equivalence, we complete the proof.
         \end{proof}

%


The following is our main theorem,
 which follows immediately from Theorems \ref{thm.hw2.1}
 and  \ref{thm.WZdinfty}. 
Let $(y_t)$ be as above. For the equal partition 
 $\cQ_m :=\{0< 1/m <\cdots < (m-1)/m<1 \}$ of length $1/m$,
 we simply write $y (m)$ for $y(\cQ_m)$ ($m \in\N$).

\begin{theorem} \label{thm.main}
Let $H \in (1/4,1/2]$, $\kappa \in (0, 2H -\tfrac12)$ and $t \in (0,1]$.
We assume that $y_t \in \mathbf{D}_\infty (\R^e)$ is non-degenerate in the sense of Malliavin.

Then, for every $r>0$, $\beta \ge 0$, $\delta >0$ and $1<q< \infty)$
  such that $r > \beta + e/q' +1$ (with $1/q +1/q' =1$), it holds that
\begin{equation}\label{eq.huwat2}
\sup_{\xi\in\R^e}  \bigl\|  
  [(1 -\Delta)^{\beta/2}\varphi_{m^{-\delta}}] (y(m)_t - \xi)
      - [(1 -\Delta)^{\beta/2}\delta_\xi] (y_t)
       \bigr\|_{ \mathbf{D}_{q, -r}} = 
                    O (m^{-\kappa \wedge \delta})
\end{equation}  
as $m\to\infty$. Recall that $\varphi_\rho (x)$ is defined by \eqref{def.HK}.
\end{theorem}

\begin{remark} 
We make a few comments on Theorem \ref{thm.main} above.
\\
{\rm (i)}~By taking a generalized expectation (i.e. the pairing 
with the constant function $\mathbf{1}$) of the left hand side of \eqref{eq.huwat2} when $\beta =0$, we obtain an approximation 
of the density of the law of $y_t$ with respect to the Lebesgue measure.
(See a comment after Theorem \ref{thm.hw2.1}, too.)
\\
{\rm (ii)}~
A typical sufficient condition for the non-degeneracy of $y_t$
is H\"ormander's bracket generating condition on 
$\{ V_0, V_1, \ldots, V_d\}$ 
at the starting point $a$. 
Here, we write $b=V_0$ and $\sigma = [V_1, \ldots, V_d]$
($V_i$'s are viewed as vector fields on $\R^e$ in this remark).
When $H=1/2$, this fact is classical.
When $H \in (1/4, 1/2)$, see \cite{chlt} (and also \cite{got2}).

The precise statement of the above condition is as follows:
Set $\Lambda_1 := \{V_i \mid 1\le i \le d\}$ and 
$\Lambda_k := \{[V_i, Z] \mid 0\le i \le d,\, Z\in \Lambda_{k-1} \}$ for $k\ge 2$.
We say that $\{ V_0, V_1, \ldots, V_d\}$ satisfies 
H\"ormander's bracket generating condition at $a$ if 
$\{ Z (a) \mid Z \in \cup_{k \in \N} \Lambda_k\}$ linearly spans $\R^e$.
(This is a traditional condition and there are many examples
which satisfy it.)
\\
{\rm (iii)}~Suppose that $\{ y(m)_t\}_m$ are  uniformly non-degenerate in the sense of Malliavin. Then, by \cite[Corollary 2.2]{hw},
$ [(1 -\Delta)^{\beta/2}\varphi_{m^{-\delta}}] (y(m)_t  -\xi)$
and $O (m^{-\kappa \wedge \delta})$ in \eqref{eq.huwat2}
can be replaced by $[(1 -\Delta)^{\beta/2}\delta_\xi] (y(m)_t)$
and $O (m^{-\kappa})$, respectively. 
However, even when $H=1/2$,
 the author does not know a suitable sufficient 
condition for the non-degeneracy of $y(m)_t$. 
\end{remark}


\noindent
{\bf Acknowledgement:}~
The author thanks Professors Arturo Kohatsu-Higa,
Nobuaki Naganuma and Setsuo Taniguchi
for their helpful comments.
He is supported by 
JSPS KAKENHI (Grant No. 20H01807).


\bigskip
\begin{flushleft}
  \begin{tabular}{ll}
    Yuzuru \textsc{Inahama}
    \\
    Faculty of Mathematics,
    \\
    Kyushu University,
    \\
    744 Motooka, Nishi-ku, Fukuoka, 819-0395, JAPAN.
    \\
    Email: {\tt inahama@math.kyushu-u.ac.jp}
  \end{tabular}
\end{flushleft}


\begin{thebibliography}{00}

\bibitem{br}
Bally, V.; Rey, C.; 
Approximation of Markov semigroups in total variation distance. 
Electron. J. Probab. 21 (2016), Paper No. 12, 44 pp. 

\bibitem{bt}
Bally, V.; Talay, D.; 
The law of the Euler scheme for stochastic differential equations. II. Convergence rate of the density. 
Monte Carlo Methods Appl. 2 (1996), no. 2, 93--128. 


\bibitem{bfrs}
Bayer, C.; Friz, P. K.; Riedel, S.; Schoenmakers, J.;
From rough path estimates to multilevel Monte Carlo.
SIAM J. Numer. Anal. 54 (2016), no. 3, 1449--1483.


\bibitem{chlt}
Cass, T.; Hairer, M.; Litterer, C.; Tindel, S.;
Smoothness of the density for solutions to Gaussian rough differential equations. 
Ann. Probab. 43 (2015), no. 1, 188--239. 

\bibitem{cll}
Cass, T.; Litterer, C.; Lyons, T.; 
Integrability and tail estimates for Gaussian rough differential equations.
Ann. Probab. 41 (2013), no. 4, 3026--3050. 

\bibitem{fggr}
Friz, P.; Gess, B.; Gulisashvili, A.; Riedel, S.;
The Jain-Monrad criterion for rough paths and applications to random Fourier series and non-Markovian H\"ormander theory. 
Ann. Probab. 44 (2016), no. 1, 684--738. 


%
\bibitem{fr}
Friz, P.;  Riedel, S.;
Convergence rates for the full Gaussian rough paths.
Ann. Inst. Henri Poincar\'e Probab. Stat. 50 (2014), no. 1, 154--194. 

\bibitem{fvbk}
Friz, P.; Victoir, N.; {\it Multidimensional stochastic processes as rough paths.} Cambridge University Press, Cambridge, 2010.

\bibitem{got2}
Geng, X.; Ouyang, C.; Tindel, S.;
Precise local estimates for differential equations driven by fractional Brownian motion: hypoelliptic case. 
Ann. Probab. 50 (2022), no. 2, 649--687. 

\bibitem{hu}
Hu, Y.;
{\it Analysis on Gaussian spaces.} World Scientific, 2017.


\bibitem{hw}
Hu, Y.; Watanabe, S.;
Donsker's delta functions and approximation of heat kernels by the time discretization methods. 
J. Math. Kyoto Univ. 36 (1996), no. 3, 499--518. 

\bibitem{iwbk}
Ikeda, N., Watanabe, S.; {\it Stochastic differential equations and diffusion processes. }
Second edition. 
North-Holland Publishing Co., Amsterdam; Kodansha, Ltd., Tokyo, 1989.

\bibitem{ina14}
Inahama, Y.;  
Malliavin differentiability of solutions of rough differential equations. 
J. Funct. Anal. 267 (2014), no. 5, 1566--1584. 


\bibitem{kh}
Kohatsu-Higa, A.;
High order It\^o-Taylor approximations to heat kernels. 
J. Math. Kyoto Univ. 37 (1997), no. 1, 129--150. 

\bibitem{kh01}
Kohatsu-Higa, A.;
Weak approximations. A Malliavin calculus approach. 
Math. Comp. 70 (2001), no. 233, 135--172. 

\bibitem{lcl}
 Lyons, T.; Caruana, M.; L\'evy, T.; 
{\it Differential equations driven by rough paths. }
  Lecture Notes in Math., 1908. Springer, Berlin, 2007.

\bibitem{lq}
 Lyons, T.; Qian, Z.; 
{\it System control and rough paths. }
Oxford University Press, Oxford, 2002.

\bibitem{mt}
Matsumoto, H.; Taniguchi, S.;
{\it Stochastic analysis. 
It\^o and Malliavin calculus in tandem. }
Cambridge University Press, Cambridge, 2017.

\bibitem{na}
Naganuma, N.;
Exact convergence rate of the Wong-Zakai approximation to RDEs driven by Gaussian rough paths. 
Stochastics 88 (2016), no. 7, 1041--1059. 

\bibitem{nu}
 Nualart, D.; 
 {\it The Malliavin calculus and related topics. }
 Second edition. Springer-Verlag, Berlin, 2006. 



\bibitem{sh}
Shigekawa, I.; {\it Stochastic analysis.}
Translations of Mathematical Monographs, 224. Iwanami Series in Modern Mathematics. 
American Mathematical Society, Providence, RI, 2004. 


\bibitem{str}
Stroock, D. W.;
{\it Gaussian Measures in Finite and Infinite Dimensions.}
Springer Verlag, Cham, 2023.

\end{thebibliography}
\end{document}